\documentclass[preprint,number]{elsarticle}

 \usepackage{amsmath}
 \usepackage{amsthm}
 \usepackage{amssymb}
 \usepackage{eucal}

 \usepackage[all]{xy}
 \usepackage{xypic}
 \theoremstyle{plain}
 \newtheorem{theorem}{Theorem}[section]
 \newtheorem*{maintheorem}{Main Theorem}
 \newtheorem{proposition}[theorem]{Proposition}
 \newtheorem{lemma}[theorem]{Lemma}
 \newtheorem{corollary}[theorem]{Corollary}

 \theoremstyle{definition}
 \newtheorem{definition}[theorem]{Definition}
 \newtheorem{convention}[theorem]{Convention}
\newtheorem{remark}[theorem]{Remark}
 \newtheorem{ex}[theorem]{Example}
 \newtheorem{notation}[theorem]{Notation}

 \newcommand{\nc}{\newcommand}
 \nc{\op}{\operatorname}

 \nc{\mb}{\mathbb}
 \nc{\mc}{\mathcal}
 \nc{\mf}{\mathfrak}

 \nc{\Hom}{\op{Hom}}
 \nc{\Ext}{\op{Ext}}

 \nc{\Sym}{\op{Sym}}

 \nc{\Aut}{\op{Aut}}

 \nc{\Spec}{\op{Spec}}
 \nc{\Proj}{\op{Proj}}
 \nc{\Pn}{\mathbb{P}}

 \nc{\Fun}{\op{Fun}}
 \nc{\mO}{\mc{O}}

 \nc{\Z}{\mb{Z}}
 \nc{\N}{\mb{N}}
 \nc{\Q}{\mb{Q}}
 \nc{\R}{\mb{R}}
 \nc{\C}{\mb{C}}
 \nc{\Cs}{\mb{C}^{\times}}

 \nc{\pt}{\op{pt}}

 \nc{\Mbar}{\overline{\mc{M}}}
 \nc{\Mtwid}{\widetilde{\mc{M}}}

 \nc{\Bun}{\op{Bun}}
 \nc{\Pic}{\op{Pic}}
 \nc{\Jac}{\op{Jac}}

 \nc{\mP}{\mc{P}}
 \nc{\wmP}{\widetilde{\mP}}
 \nc{\XP}{X_{\mP}}
 \nc{\mL}{\mc{L}}
 \nc{\md}{d} 

 \nc{\un}{n}
 \nc{\dt}{\mathfrak{t}}
 \nc{\uu}{\underline{u}}
 \nc{\ul}{\underline{l}} 
 \nc{\uN}{\underline{N}}
 \nc{\uM}{\underline{M}}
 \nc{\ev}{\op{ev}}
 \nc{\st}{\op{st}}
 \nc{\id}{\op{id}}
 \nc{\dct}{\op{def}}

\begin{document}

\title{Gromov-Witten Gauge Theory}

\address[berkeley]{Department of Mathematics, University of California, Berkeley, CA 94720, USA}
\address[stonybrook]{Simons Center for Geometry \& Physics, Stony Brook University, Stony Brook, NY, 11794}

\author[berkeley]{Edward Frenkel\fnref{fn1}}
\ead{frenkel@math.berkeley.edu}

\author[berkeley]{Constantin Teleman\fnref{fn1}}
\ead{teleman@math.berkeley.edu}

\author[stonybrook]{A.J. Tolland\fnref{fn1}}
\ead{ajt@math.sunysb.edu}

 \fntext[fn1]{Partially supported by DARPA and AFOSR through the grant
 FA9550-07-1-0543}

 \begin{abstract}
 We introduce a modular completion of the stack of maps from stable
 marked curves to the quotient stack $[\pt/\Cs]$, and use this stack to construct
 some gauge-theoretic analogues of the Gromov-Witten invariants.  
 We also indicate the generalization of these invariants to the quotient stacks $[X/\Cs]$, where $X$ is a smooth proper complex algebraic variety.
 \end{abstract}

 \maketitle 

In this paper we construct algebro-geometric Gromov-Witten invariants for the quotient stack $[\pt/\Cs]$.  These invariants are the ``twisted'' indices of products of evaluation K-theory classes on certain moduli stacks $\Mtwid_{g,I}([\pt/\Cs])$, which classify maps from marked nodal curves to the quotient stack $[\pt/\Cs]$.  These are the moduli stacks of principal $\Cs$-bundles on such curves.  We do not impose stability conditions on these bundles, and our moduli stacks are Artin stacks that are not  proper.  Consequently, the existence of these invariants is non-trivial.  Our main theorem asserts that they are, in fact, well-defined.  

 This construction is the first step in a larger project, the goal of which is to define Gromov-Witten invariants for the Artin stacks $[X/G]$, where $X$ is a smooth projective variety and $G$ is a reductive algebraic group over $\C$.  Similar invariants (in the case where $G$ is a finite group) have already been defined in the symplectic setting by Chen \& Ruan \cite{MR1950941} and algebro-geometrically by Abramovich, Graber, Olsson, \& Vistoli \cite{MR1862797, MR2450211}.  In this paper, we explain (modulo a certain conjecture) how to define invariants for $[X/\Cs]$ in terms of our invariants for $[\pt/\Cs]$.

 We expect that these invariants -- defined as twisted indices of products of evaluation classes on moduli stacks of maps from marked curves to $[X/\Cs]$ -- may be interpreted as correlation functions in a topological quantum field theory.  This quantum field theory is a {\it gauge} theory.  The stack of algebraic principal $G$-bundles on a smooth curve $\Sigma$ is homotopy equivalent \cite{MR702806} to the stack of $K$-connections on $\Sigma$, where $K$ is the compact form of $G$, and so one can view our construction of Gromov-Witten invariants for $[\pt/\Cs]$ as a topological version of Feynman path integration over the space of $U(1)$-connections.

 Similar invariants have been studied by Mundet, Cieliebak, Gaio, Salamon, and Tian 
 \cite{MR1801657, MR1777853,  MR1959059, MR1953239, Riera:math0404407}, who defined 
 invariants for symplectic manifolds with Hamiltonian $U(1)$-actions by integrating 
 over coarse moduli spaces of solutions to the symplectic vortex equations, and by 
 Gonzalez \& Woodward \cite{arXiv:0811.3358, Gonzalez:arXiv0907.3869}, who considered 
 the generalization from $U(1)$ to compact reductive groups.  In fact, the idea that 
 the path integral of a gauged nonlinear sigma model might compute cohomological invariants 
 of moduli spaces of connections $A$ and $\overline{\partial}_A$-holomorphic sections of 
 a bundle with fiber $X$ appears in Witten's original paper on topological sigma models 
 \cite{Witten:1988xj}.  For a more recent perspective on this, see the study \cite
 {frenkel-2008} of gauged sigma models in the ``infinite radius limit''. 

 Our invariants differ from the existing ones in that we do not impose stability conditions in the sense of geometric invariant theory on the maps to $[X/\Cs]$ and the invariants take values in $K$-theory rather than cohomology.  (Our approach has been outlined, in particular, in a lecture given by one of us at the Seattle conference in 2005 \cite{TelAG}.)  The proof of finiteness of these invariants uses a technical argument adapted from \cite{MR1792291}.  Consequently, our invariants are defined in greater generality than the known ones.  Indeed, we expect to recover the Gromov-Witten invariants of GIT quotients from our invariants by applying the Chern character to certain limits of our invariants.  This was done for smooth curves and 
 $G$-bundles in \cite{math.AG/0312154}.

 \subsection{Sketch of the Construction}  We explain here the construction of invariants 
 for $[\pt/\Cs]$. 

Recall that a $\Cs$-bundle on a nodal curve $\Sigma$ is defined by a $\Cs$-bundle on 
the normalization of $\Sigma$ together with an identification of the two fibers at the 
preimages of each node. The stack $\Bun_{\Cs}(g,I)$ of $\Cs$-bundles over the universal 
stable curve fails to be complete, because the space of identifications over a given node 
is isomorphic to $\Cs$.  Following Gieseker \cite{MR739786} and Caporaso \cite{MR1254134}, 
we add new strata which represent the limits where an identification goes to zero or infinity, 
by allowing projective lines carrying the line bundle ${\mathcal O}_{\mb{P}^1}(1)$ to appear 
at the nodes. The resulting stack -- denoted $\Mtwid_{g,I}([\pt/\Cs])$ -- is complete but 
not  separated, i.e., the limit of a family of bundles exists but may not be unique.  
(Similar completions of various stacks of vector bundles on nodal curves have
 been studied by several authors, see \cite{MR771150,
   MR2105707,MR2735773}.  Our definition was inspired by Caporaso's
 thesis \cite{MR1254134} and the papers of Nagaraj \& Seshadri
 \cite{MR1455315, MR1687729}.)


The stack $\Mtwid_{g,I}([\pt/\Cs])$ has a forgetful map
\[
F: \Mtwid_{g,I}([\pt/\Cs]) \to \Mbar_{g,I}.
\]
It also carries a universal curve $\pi: C \to \Mtwid_{g,I}([\pt/\Cs])$ with marked sections 
$\sigma_i$ and a universal principal $\Cs$-bundle $p: \mP \to C$. A $\Cs$-representation $V$
leads to an associated vector bundle $\phi^*V$ on $C$ (pulled back by the classifying map 
$\phi: C \to [\pt/\Cs]$ of $\mP$, in the language of stacks).  Restricting  to the $\sigma_i$, 
we obtain the {\it evaluation classes} $\ev_i^*[V]$ in the K-theory of $\Mtwid_{g,I}([\pt/\Cs])$.
 
The Gromov-Witten invariants of $[\pt/\Cs]$, like those of a variety $X$, result from pushing 
a product of evaluation classes forward along the forgetful morphism $F$.\footnote{The stack 
$\Mtwid_{g,I}([\pt/\Cs])$ turns out to be unobstructed, so virtual structure theory is not 
required.} However, our setup differs from the standard one in two ways.
 \begin{enumerate}
 \item Our invariants are constructed in K-theory, rather than cohomology.
 \item Our invariants are always {\it twisted}, in the sense of \cite{MR2276766}.
 \end{enumerate}

 Twisting requires some definition.  A line bundle $\mL$ on $\Mbar_{g,I}([\pt/\Cs])$ is 
 {\it admissible} if 
 \[
 \mL \simeq (\op{det}R\pi_*\phi^*\C_1)^{\otimes (-q)},
 \]
where $\C_1$ is the standard representation and $q$ is a positive rational number.  
An {\it $\mL$-twisted Gromov-Witten invariant of $[\pt/\Cs]$} is the K-theory class in 
$\Mbar_{g,I}$ of the total direct image along $F$ of the twist by (admissible) $\mL$ of 
a product of evaluation classes:
\[
RF_* \big(\mL \bigotimes \otimes_{i \in I} \ev_i^*V_i \big)  
\]
 One can generalise this by inserting products of {\it index classes $R\pi_*\phi^*V$ along $C$}; 
these classes may be assembled into {\it higher twistings} \cite{MR2276766}.  
One also defines gravitational descendants, tensoring each evaluation class with a 
power of the tangent line $T_i$ to $C$ at $\sigma_i$. An {\it admissible complex} is a 
sum of complexes of the form 
 \[
 \mL \bigotimes \otimes_a (R\pi_*\phi^*_{a}) \bigotimes \otimes_i (\ev_i^*V_{i} \otimes T_i^{\otimes n_i}).
 \]
The subring of $K(\Mtwid_{g,I}([\pt/\Cs])$ generated by such products is called the {\it ring 
of admissible classes}.  It is a subring without unit; the trivial line bundle $\mc{O}$ is not admissible.   

It is not clear that the push-forward of an admissible class along $F$ is well-defined, because the moduli stacks $\Mtwid_{g,I}([\pt/\Cs])$ differ from Kontsevich's stack of stable maps in two important ways.

First, $\Mtwid_{g,I}([\pt/\Cs])$ is an Artin stack, rather than a Deligne-Mumford stack; points 
can have continuous stabilizers.  This makes it implausible to integrate cohomology classes 
along the morphism $F$, and it is for this reason that we use K-theory instead of cohomology.  
The problem 
can be understood as follows. The fibers of $F$ are quotients stacks of the form $[A/\mc{G}]$ 
($\mc{G}$ a group), and integrating over $[A/\mc{G}]$  factors into the steps 
$[A/\mc{G}] \to [\pt/\mc{G}]\to \pt$. The first step may be sensible, but the second integration 
gives zero in cohomology, by reason of degree: it shifts the degree of a cohomology class by 
$-\dim([\pt/\mc{G}]) = \dim(\mc{G})$, so the only classes that could survive come from 
the zero group $H^{-\dim(\mc{G})}_\mc{G}(\pt)$.  
In K-theory, on the other hand, the pushforward along $[\pt/\mc{G}] \to \pt$ does exist: 
for reductive $\mc{G}$, it is implemented by the functor of 
$\mc{G}$-invariants, from $\mc{G}$-representations to vector spaces.

The second problem is that $\Mtwid_{g,I}([\pt/\Cs])$ is not proper.  It is complete, but it is in general neither separated nor of finite type.   Thus, the existence of a pushforward along the forgetful morphism to $\Mbar_{g,I}$ is a delicate matter; not every K-theory class on $\Mbar_{g,I}([\pt/\Cs])$ has a well-defined index.  The main theorem in this paper asserts that the index map is well-defined for the admissible classes.
 
 \begin{maintheorem}\label{maintheorem}
 The derived pushforward $RF_*\alpha$ of an admissible complex $\alpha$ along the bundle-forgetting 
 map $F: \Mtwid_{g,I}([\pt/\Cs]) \to \Mbar_{g,I}$ is a (bounded) complex of coherent sheaves.
 \end{maintheorem}
\noindent
 This theorem is a relative version (over varying curves) of the analogous finiteness result for 
 $\Bun_G(\Sigma)$ in 
 \cite{math.AG/0312154}.  The proof, in rough outline:  
 \begin{enumerate}
 \item The restriction  of $\Mtwid_{g,I}([\pt/\Cs])$ to a small enough affine \'etale chart $B$ of 
 $\Mbar_{g,I}$ 
 will be presented as a stack quotient  $[A/\mc{G}]$.  Here, $A$ is an algebraic space, the moduli 
 space of Gieseker bundles trivialized  at certain special points, and $\mc{G}$ is the group of 
 rescalings of these trivializations.

\item  We identify an open  subspace $A^o \subset A$ for which $[A^o/\mc{G}] \simeq Q 
\times [\pt/\Cs]$, with $Q$ proper over $B$. This pertains to a $\mc{G}$-equivariant 
stratification of  $A$ which will let us relate cohomologies over $A^o$ and $A$.  
 
\item  By estimating the $\mc{G}$-weights on the fixed point fibers of admissible classes, 
we show that the $\mc{G}$-invariants in the local cohomologies of $\alpha$ on the strata 
above are always finite over $B$, and almost always vanish.  This shows 
the finiteness over $B$ of the $\mc{G}$-invariants  in the global sections $R\Gamma(A,\alpha)$ 
and proves the theorem, since $\mc{G}$ is reductive.  
 \end{enumerate}

 \subsection{Invariants for $[X/\Cs]$}  
The stack $[X/\Cs]$ is defined so that maps from a curve $\Sigma$ to $[X/\Cs]$ correspond to pairs 
$(\mP,s)$ consisting of a principal $\Cs$-bundle and a section $s \in \Gamma(\Sigma,\mP \times_{\Cs} X)$ of the associated  bundle with fiber $X$.  To define Gromov-Witten invariants for $[X/\Cs]$,  
we need a moduli stack $\Mtwid_{g,I,\beta}([X/\Cs])$ of curves and degree $\beta$ maps to 
$[X/\Cs]$  which supports the tautological classes, and a push-forward operation on these classes 
from  $\Mtwid_{g,I,\beta}([X/\Cs])$ to $\Mbar_{g,I}$.  

In the final section, we define a moduli stack $\Mtwid_{g,I,\beta}([X/\Cs])$, along with a natural 
section-forgetting morphism
\[
F_\beta: \Mtwid_{g,I,\beta}([X/\Cs]) \to \Mtwid_{g,I}([\pt/\Cs]),
\]
whose fibers are stacks of sections of bundles with fiber $X$ (associated to Gieseker 
bundles).  Sections are locally maps and can develop local singularities in the same way. 
Following Kontsevich, we ensure that $F_\beta$ is proper by allowing bubbling at 
singularities. Moreover, we also show that $F_\beta$ is Deligne-Mumford 
and carries a perfect obstruction theory. This makes our morphism is very much like the 
$F_\beta: \Mbar_{g,I,\beta}(X) \to \Mbar_{g,I}$ of ordinary Gromov-Witten theory 
(and our proofs are straightforward variations of the standard ones).  

These facts imply the existence of a virtual K-theoretic pushforward along $F_\beta$.  
We expect that the image of an admissible class can be further integrated down to 
$\Mbar_{g,I}$ --- namely, that it satisfies  bounds for the weights at fixed-points 
similar to those we check for the admissible classes on $\Mtwid_{g,I}([\pt/\Cs])$. 
If so, this would result 
in Gromov-Witten invariants of $[X/\Cs]$ in the K-theory of $\Mbar_{g,I}$. 

 \subsection{Plan of the Paper}

Section \ref{thestack}  reviews basic facts about nodal 
curves and principal $\Cs$-bundles. The moduli stack $\Mtwid_{g,I}([\pt/\Cs])$ of 
Gieseker bundles on stable curves is introduced with some key examples (small $g$ and $|I|$).

In Section \ref{ElemGeom}, we prove some basic facts about the geometry of our stack: 
it is an Artin stack, is stratified by topological type, and is complete 
(but not separated).

In Section \ref{usefulatlas}, we give an (\'etale-local) presentation of 
$\Mtwid_{g,I}([\pt/\Cs])$  as a quotient $A/\mc{G}$ (where $\mc{G} \simeq (\Cs)^V$). 
We identify a stable subspace $A^o\subset A$ which leads to a smooth and proper 
quotient moduli space over  $\Mbar_{g,I}$.

In Section \ref{fixedpoints}, we refine the stratification by topological type by 
tracking the nodes smoothed under deformations. We use this to stratify $A/\mc{G}$ 
by distinguished spaces $Z,W$ which are affine space bundles over their fixed-point 
loci under subgroups of $\mc{G}$. 

In Section \ref{admissibleclasses}, we review the admissible K-theory classes and estimate 
the weights of the fixed-point fibers of subgroups of $\mc{G}$.  

In Section \ref{theinvariants},  we use a local cohomology vanishing 
argument to we finish the proof of the main theorem.

 In Section \ref{XmodC}, we construct a moduli stack which we expect to carry Gromov-Witten invariants 
 for $[X/\Cs]$.

 \subsection{Acknowledgements}

E.F. thanks Andrei Losev and Nikita Nekrasov for useful discussions.  
A.T. thanks Jarod Alper, Tom Coates, Dan Edidin, Ezra Getzler, Eduardo Gonzalez, Reimundo Heluani, 
and Chris Woodward for helpful conversations, and also the referee for a nice hairdryer treatment 
concerning the exposition.

 This research has been partially supported by DARPA and AFOSR through the grant FA9550-07-1-0543. In addition, E.F. and A.T. were supported by the NSF grant DMS-0303529, and C.T. was supported by the NSF grant DMS-0709448.

 \section{The Stack of Gieseker Bundles}\label{thestack}
In this section, we introduce appropriate moduli stacks of marked curves carrying principal 
$\Cs$-bundles and discuss a few simple examples.

 \subsection{Curves}
We always work over $\C$. In everything that follows, $(C,\sigma_i)$ is a family of prestable marked curves over a finitely-generated complex base scheme $B$.  More precisely, $\pi: C \to B$ is a flat proper morphism whose fibers are connected complex projective curves of genus $g$ with at worst nodal singularities, carrying a collection of smooth marked points $\sigma_i: B \to C$ which are indexed by an ordered set $I$.  A point is 
{\it special} if it is a node or a marked point.  Special points are required to be pairwise disjoint. 
We shall always assume that any rational component of $C$ has at least {\it two} special points.

We  reserve the notation $(\Sigma,\sigma_i)$ for families of {\it stable} marked curves.  
Recall that a marked curve is stable if each component of genus 0 carries at least 3 special points and each component of genus 1 carries at least 1 special point.  
The {\it stabilization morphism} $\st: C \to C^{st}$ blows down every unstable rational curve in $C$.  Stabilization can be implemented by a pluricanonical embedding and thus works in families. 

In addition to marked points and nodes, we will use a third kind of special points on $C$:  
{\it trivialization points}.  These are smooth points at which we will trivialize the fibers of 
a $\Cs$-bundle, but which do not count in determining stability. No harm comes if stabilization 
points agree with marked points or with each other, but none results from ruling this out, either. 

\begin{notation}
We denote all special points by $\sigma$, distinguishing them by the subscript. Ordinary 
marked points are denoted $\sigma_i$, with $i\in I$.   Nodes are $\sigma_e$, with $e$ in a set 
$E$.  Trivialization points are denoted $\sigma_v$, with $v$ in a set $V$. 
 \end{notation}

\begin{figure}[htb] \label{graphpict}
\[\begin{xy}
(-20,0)*{\begin{xy}
(0,0)*\xycircle(10,7){-};
(6,0)*\xycircle(4,2){-};
(-3,4)*{\bullet};
(-1,2)*{\sigma};
(-8,4);(-8,-4) **\crv{(-12,6) & (-15,0) & (-12,-6)};
(-8,4);(-8,-4) **\crv{(-10,8) & (-20,15) & (-40,0) & (-20,-15)& (-10,-8)};
(-25,0)*{\begin{xy}
(5,1);(-5,1) **\crv{(0,-4)};
(4,0);(-4,0) **\crv{(0,2)};
\end{xy}};
\end{xy}};
(30,0)*{\begin{xy}
(0,0)*\xycircle(4,4){-} *{1};
(20,0)*\xycircle(4,4){-} *{0};
(2.83,2.83); (17.17,2.83)**\crv{(10,7)};
(2.83,-2.83); (17.17,-2.83)**\crv{(10,-7)};
(20,5);(20,10)**\crv{};
(22.83,2.83);(22.83,-2.83) **\crv{(30,6) & (35,0)& (30,-6)};
\end{xy}};
\end{xy}\]
\caption{A (stable) marked curve and its associated modular graph, which has two splitting edges, 
one self-edge and one tail.}
\end{figure}
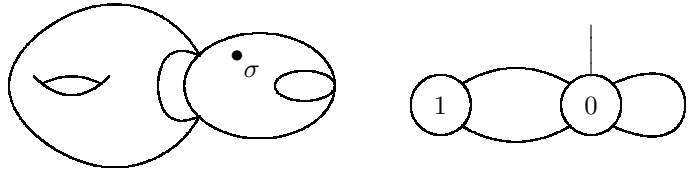

The topology of a marked curve can be encoded in its {\it modular graph} \cite{MR1412436}. 
This consists of a finite, undirected graph $\gamma$, with vertex set $V_\gamma$ and edge set 
$E_\gamma$, together with a function $g:V_\gamma\to \N$. There is one vertex for each component, 
one full edge for each node and one half-edge (or {\it tail}) for each special point; 
$g(v)$ is the arithmetic genus of the irreducible component labeled by $v$. Multiple edges 
and loops ({\it self-edges}) are permitted, matching the geometry of the curve, as 
in Figure~\ref{graphpict}. The underlying topological space of 
$\gamma$ is denoted $|\gamma|$.

Full subgraphs of $\gamma$ correspond to unions of components of the curve; leaving out 
full edges leads to partial normalizations. A subgraph normally carries the relevant old tails, 
along with new ones at the severed edges. Attributes of the curve transfer to the modular 
graph without comment; we thus have stable graphs and subgraph, an (arithmetic) genus of a 
graph, and self- or splitting edges according to the type of the node. 

$\Mbar_{g,I}$ denotes the stack of stable genus $g$, $I$-marked curves. Its locally closed 
substack which classifies marked curves of type $\gamma$ is denoted $\mc{M}_\gamma$.  
These substacks stratify $\Mbar_{g,I}$ in a strong sense: the boundary of any $\mc{M}_\gamma$ 
is a union of strata $\mc{M}_{\gamma'}$.  The boundary divisors in this stratification 
have normal crossings.  Modular graphs labelling strata are always connected.  

\begin{convention}
In this paper, {\it strata} are locally closed, 
but generally not closed.  $\mc{M}_\gamma$ is exactly the stack of curves of topological type 
$\gamma$, not its closure.
\end{convention}

\subsection{Bundles on Nodal Curves}\label{bunnodecurve}

\begin{definition}
A {\it principal $\Cs$-bundle on a scheme $X$} is a scheme $\mP$ on which $\Cs$ acts 
freely (from the right) and a $\Cs$-invariant map $p: \mP \to X$ which is locally\footnote
{Here, we may use the Zariski topology, but \'etale covers will be more common in 
the paper.} trivial: $X$ has an open cover $\{U_\alpha\}$ such that $U_\alpha \times_p 
\mP \simeq U_\alpha \times \Cs$, $\Cs$-equivariantly. 
A {\it family of principal $\Cs$-bundles on $X$} (parametrized by a scheme $B$) is a 
principal bundle on $X \times B$.
\end{definition}

\begin{definition}
The {\it degree} of a principal $\Cs$-bundle $\mP$ over an irreducible curve $X$ 
is the Chern class, in $H^2(X;\Z)\simeq \Z$, of the associated line bundle. 

For a pre-stable curve $C$ with modular graph $\gamma$, the {\it multi-degree} $\md$ of 
$\mP$ assigns to each irreducible component $C_v \subset C$ the degree $d_v$ of $\mP|_{C_v}$:     
\[
\md: V_\gamma \to \Z, \qquad v  \mapsto d_v.
\]
The {\it total degree $D$} is the sum $D = \sum_{v\in V} d_v$ .  The {\it topological type} of 
$\mP\to C$ is the pair $(\gamma,\md)$.
\end{definition}

\begin{remark}\label{normalize}
Let $\nu: \widetilde{C} \to C$ denote the normalization of a nodal curve $C$.  
A principal bundle $\mP$ on $C$ is equivalent to the following data:
\begin{enumerate}
\item a principal $\Cs$-bundle $\wmP$ ($=\nu^*\mP$) on $\widetilde{C}$, and
\item for each node $\sigma \in \Sigma$, a {\it gluing isomorphism} $\iota: \wmP_{\sigma^+} \simeq \wmP_{\sigma^-}$, which identifies the fibers of $\wmP$ over the preimages $\sigma^\pm$ of $\sigma$ under $\nu$.
\end{enumerate}
\end{remark}
\begin{figure}[htb]
\centering
\[\begin{xy}
(-2,0)*{\begin{xy}
(0,0)*\xycircle(12,8){-};
(-5,1)*{};(1,1)*{} **\crv{(-2,-1)};
(-4,.5)*{};(0,.5)*{} **\crv{(-2,2)};
\end{xy}};
(-2,20)*{\begin{xy}
(0,0)*{}; (24,0)*{} **\crv{(10,2)};
(0,20)*{}; (24,20)*{} **\crv{(10,22)};
(0,0)*{};(0,20)*{} **\dir{-};
(24,0)*{};(24,20)*{} **\dir{-};
\end{xy}};
(20,0)*{\begin{xy}
(0,0)*\xycircle(10,8){-};
(-4,0)*{}; (4,0) *{} **\crv{(0,-3)};
(-3,-1)*{}; (3,-1)*{} **\crv{(0,2)};
\end{xy}};
(20,20)*{\begin{xy}
(0,0)*{}; (20,0)*{} **\crv{(10,5)};
(0,20)*{}; (20,20)*{} **\crv{(10,25)};
(0,0)*{};(0,20)*{} **\dir{-};
(20,0)*{};(20,20)*{} **\dir{-};
\end{xy}};
(40,10)*{\simeq};
(60,0)*{\begin{xy}
(0,0)*\xycircle(12,8){-};
(-5,1)*{};(1,1)*{} **\crv{(-2,-1)};
(-4,.5)*{};(0,.5)*{} **\crv{(-2,2)};
(12,0)*{\bullet};
\end{xy}};
(60,20)*{\begin{xy}
(0,0)*{}; (24,0)*{} **\crv{(10,2)};
(0,20)*{}; (24,20)*{} **\crv{(10,22)};
(0,0)*{};(0,20)*{} **\dir{-};
(24,0)*{};(24,20)*{} **\dir{-};
\end{xy}};
(100,0)*{\begin{xy}
(0,0)*\xycircle(10,8){-};
(-4,0)*{}; (4,0) *{} **\crv{(0,-3)};
(-3,-1)*{}; (3,-1)*{} **\crv{(0,2)};
(-10,0)*{\bullet};
\end{xy}};
(100,20)*{\begin{xy}
(0,0)*{}; (20,0)*{} **\crv{(10,5)};
(0,20)*{}; (20,20)*{} **\crv{(10,25)};
(0,0)*{};(0,20)*{} **\dir{-};
(20,0)*{};(20,20)*{} **\dir{-};
\end{xy}};
{\ar@{->} (89,0)*{}; (73,0)*{}};
{\ar@{->} (89,20)*{}; (73,20)*{}};
(81,25)*{\iota};
(76,3)*{\sigma^+}; (88,3)*{\sigma^-};
\end{xy}\]
\caption{Realizing $\mP$ as $\widetilde{\mP}$ together with a gluing isomorphism $\iota$}
\end{figure}
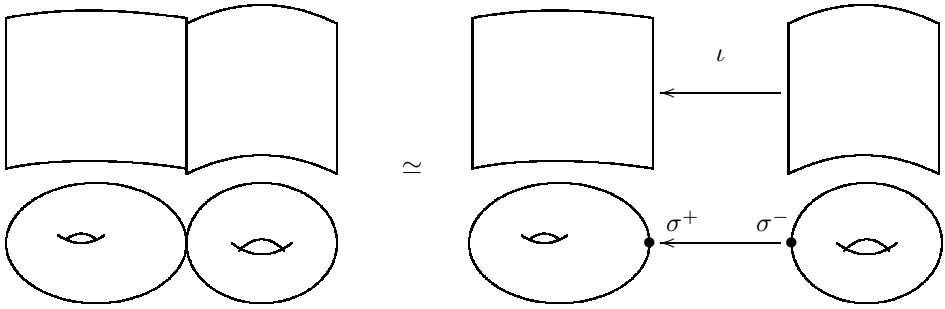

One sees from this how $\Cs$-bundles on $C$ can become singular in families: the space of gluing 
isomorphisms at a node $\sigma \in C$ is a copy of $\Cs$; in a family, these isomorphisms can 
tend to the limit points $0,\infty$.  As a result,  the stack $\Bun_{\Cs}(g,I)$ of $\Cs$-bundles 
on stable marked curves of type $(g,I)$ fails the valuative criterion for completeness. 
This will be a problem for integration of cohomology or K-theory classes.

We shall address this problem by enlarging the classification problem 
slightly, ``filling in the holes'' in $\Bun_{\Cs}(g,I)$ while keeping a principal 
bundle: we will allow copies of $\Pn^1$ to appear at the nodes of stable curves, 
and insists that these $\Pn^1$ carry degree $1$ bundles.  

\begin{definition}\label{defmod}
A morphism $m: C \to \Sigma$ of prestable curves is a {\it modification} if:
\begin{enumerate}
\item $m$ is an isomorphism away from the preimage of the nodes of $\Sigma$, and
\item the preimage under $m$ of every node in $\Sigma$ is either a node or a $\Pn^1$ with two special points.
\end{enumerate}
A {\it modification of a family} $f: \Sigma \to B$ of marked prestable 
curves is a morphism $m: C \to \Sigma$ such that, for each geometric $b \in B$, 
$m_b: C_b \to \Sigma_{f(b)}$ is a modification. 
\end{definition}

\begin{remark}\label{remmod}
\begin{enumerate}
\item Finding modifications with desirable properties --- such as, smoothness 
of the total space $C$ --- may require us to change the base $B$; the reader can be 
entrusted to write out the defining diagram. 
\item
Modifications of marked curves do not introduce $\Pn^1$s at marked points, 
only at nodes. The marked points in a family $\Sigma$  lift 
uniquely to the modification, and will sometimes be denoted by the same symbol.
\end{enumerate}
\end{remark}

\begin{figure}[htb]
\centering
\[\begin{xy}
(5,15)*{\Sigma};
(0,0)*{\begin{xy}
(-2,0)*{\begin{xy}
(0,0)*\xycircle(12,8){-};
(-5,1)*{};(1,1)*{} **\crv{(-2,-1)};
(-4,.5)*{};(0,.5)*{} **\crv{(-2,2)};
(-2,-2)*{\bullet};
(1,-3)*{\sigma_1};
\end{xy}};
(20,0)*{\begin{xy}
(0,0)*\xycircle(10,8){-};
(-4,0)*{}; (4,0) *{} **\crv{(0,-3)};
(-3,-1)*{}; (3,-1)*{} **\crv{(0,2)};
\end{xy}};
\end{xy}};
(75,15)*{C};
(70,0)*{\begin{xy}
(-2,0)*{\begin{xy}
(0,0)*\xycircle(12,8){-};
(-5,1)*{};(1,1)*{} **\crv{(-2,-1)};
(-4,.5)*{};(0,.5)*{} **\crv{(-2,2)};
(-2,-2)*{\bullet};
(1,-3)*{\sigma_1};
\end{xy}};
(15,1)*{\begin{xy}
(0,0)*\xycircle(5,5){-};
\end{xy}};
(30,0)*{\begin{xy}
(0,0)*\xycircle(10,8){-};
(-4,0)*{}; (4,0) *{} **\crv{(0,-3)};
(-3,-1)*{}; (3,-1)*{} **\crv{(0,2)};
\end{xy}};
\end{xy}};
\end{xy}\]
\caption{A nodal curve $\Sigma$ with a marked point $\sigma_1$ and its unique non-trivial modification $C$}
\end{figure}
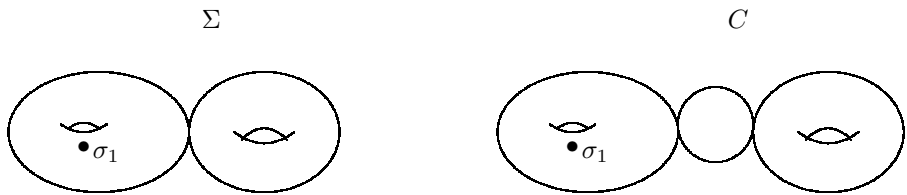

\begin{definition}[Gieseker Bundle]
Let $(\Sigma_i,\sigma_i)$ be a stable marked curve.   A {\it Gieseker $\Cs$-bundle on $(\Sigma,\sigma_i)$} is a pair $(m,\mP)$, consisting of
\begin{enumerate}
\item a modification $m: (C,\sigma_i) \to (\Sigma,\sigma_i)$  and
\item a principal $\Cs$-bundle $p: \mP \to C$,
\end{enumerate}
which satisfy the {\it Gieseker condition}:
\begin{enumerate}
\item the restriction of $\mP$ to every unstable $\Pn^1$ has degree $1$.  
\end{enumerate}
The {\it (multi)degree}  of a Gieseker bundle is the (multi)degree of the bundle $\mP$.
We will use the term {\it Gieseker bubble} for the the unstable $\Pn^1$s above.  

\medskip
If $C \to B$ is a family of prestable marked curves (with marked points $\sigma_i$), 
$\Sigma \to B$ is its stabilization, and $\mP$ a bundle on $C$ satisfying the 
Gieseker condition on all geometric fibers, we call $(C,\sigma_i,\mP)$ a 
{\it family of Gieseker bundles} on $\Sigma$.
\end{definition}

\begin{remark}
We call these Gieseker bundles after \cite{MR739786}, which was the first 
paper to describe moduli spaces parametrizing such bundles. As far as we know, 
Caporaso's [9] was the first systematic treatment of families of curves.
\end{remark}

\begin{remark} The modification is always the stabilization map of $C$, so may be 
omitted from the notation. However, different Gieseker bundles on $\Sigma$ can live 
on  different modifications $C$, so  some care is needed with the terminology.
\end{remark}

\begin{remark}\label{automph}
The Gieseker condition on a single curve $C$ can be phrased in terms of 
the topological type $(\gamma,d)$: any unstable vertices have genus $0$, originate 
two edges, and carry degree $1$. We can read a bit more from topology.

The twice marked smooth rational curve $(\Pn^1, 0, \infty)$ has a one-parameter 
family of scaling automorphisms. These can be lifted to 
any principal $\Cs$-bundle $\mP$ on $\Pn^1$; more precisely, the automorphism group 
of the marked curve with the bundle sits in the middle of an exact sequence
\[
1 \to \Cs\simeq\Aut(\mP) \to \Aut(\Pn^1,0,\infty;\mP) \to \Aut(\Pn^1,0,\infty)\simeq\Cs \to 1.
\]
The special feature of $\mP=\mc{O}(1)$ is that the automorphism group can 
be identified with $\Cs_0\times\Cs_\infty$, the group of fiber rescalings  
over $0$ and $\infty$. 

For a Gieseker bundle $\mP$ on a {\it fixed} $\Sigma$, the automorphism group of its 
restriction to any Gieseker bubble is thus identified with the $(\Cs)^2$ group of fiber 
rescalings at the nodes. No other automorphisms of $C$ are permitted, although we can 
scale the fibers of $\mP$ on the stable components. The automorphism group of $\mP$ 
is determined by $\gamma$, and agrees with the automorphism group of its restriction 
to the curve obtained from $C$ by deleting all Gieseker bubbles. 
\end{remark}

\begin{remark}\label{usefullie}
We offer some motivation for this definition, in the form of (1) an intuitive 
explanation and (2) a precise result.
\begin{enumerate}
\item Fix a nodal curve $\Sigma$ and a family of $\Cs$-bundles on $\Sigma$,  
parametrized by a coordinate $t$, for which the gluing isomorphism $\iota: \mP_{\sigma^+} \to 
\mP_{\sigma^-}$ over a node $\sigma \in \Sigma$ approaches $0$ at $t \to 0$.  (Assuming that 
$\iota \to 0$ is no loss of generality; the other limit $\iota \to \infty$ is equivalent to 
$\iota^{-1} \to 0$.)  We want to replace this singular limit 
with a bundle defined on some modification $C$. 

Consider for this this a section $s$ of the associated line bundle 
$V = \mP \times_{\Cs} \C$.  The lift $\tilde{s}$ of $s$ to the normalization must obey
\[
\tilde{s}(\sigma^+) = \iota \tilde{s}(\sigma^-).
\]
As $t \to 0$, we must have $\tilde{s}(\sigma^+) \to 0$.  By continuity, the section 
$s$ on $\Sigma$ must have a zero $z$ which approaches the node as $\iota \to 0$.  
(When $\iota \to \infty$, a zero approaches the node from the other side.)

To keep track of how a single zero $z$ approaches the node, we conformally bubble 
out a neighborhood of $\sigma^+$, creating a $\Pn^1$ at $t=0$.  The limiting 
section $s$ on this new $\Pn^1$ has one zero and no poles, so the new bundle 
on this component has degree $1$, while the degree on the original 
component will drop by $1$.

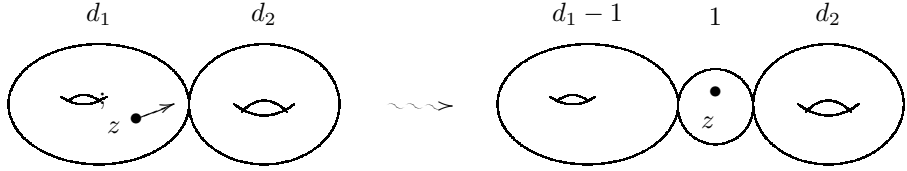
\begin{figure}[htb]
\centering
\[\begin{xy}
(0,0)*{\begin{xy}
(-2,0)*{\begin{xy}
(0,12)*{d_1};
(0,0)*\xycircle(12,8){-};
(-5,1)*{};(1,1)*{} **\crv{(-2,-1)};
(-4,.5)*{};(0,.5)*{} **\crv{(-2,2)};
(5,-2)*{\bullet};(2,-3)*{z}
{\ar (5,-2)*{};(10,0)*{}};
\end{xy}};
(20,0)*{\begin{xy}
(0,0)*\xycircle(10,8){-};
(-4,0)*{}; (4,0) *{} **\crv{(0,-3)};
(-3,-1)*{}; (3,-1)*{} **\crv{(0,2)};
(0,12)*{d_2};
\end{xy}};
\end{xy}};
{\ar@{~>} (27,-3)*{};(37,-3)*{}};
(70,0)*{\begin{xy}
(-2,0)*{\begin{xy}
(0,12)*{d_1-1};
(0,0)*\xycircle(12,8){-};
(-5,1)*{};(1,1)*{} **\crv{(-2,-1)};
(-4,.5)*{};(0,.5)*{} **\crv{(-2,2)};
\end{xy}};
(15,1)*{\begin{xy}
(0,0)*\xycircle(5,5){-};
(0,12)*{1};
(0,2)*{\bullet};(-1,-2)*{z};
\end{xy}};
(30,0)*{\begin{xy}
(0,0)*\xycircle(10,8){-};
(0,12)*{d_2};
(-4,0)*{}; (4,0) *{} **\crv{(0,-3)};
(-3,-1)*{}; (3,-1)*{} **\crv{(0,2)};
\end{xy}};
\end{xy}};
\end{xy}\]
\caption{A zero approaches a node, leading to a Gieseker bubble}
\end{figure}

\item Another explanation for the Gieseker condition stems from the relation 
with torsion-free sheaves. Recall that the latter are a  natural completion 
of the stack of curves and bundles, appearing in the Hilbert scheme closure 
of the locus of stable curves with line bundles. 

\begin{proposition}
The pushforward $m_*\mc{V}$ of the line bundle $\mc{V}$ on $C$ associated to a Gieseker bundle 
is a rank 1 torsion free sheaf on $\Sigma$.
\end{proposition}
\noindent
This fails if we allow higher degrees for $\mc{V}$ on unstable curve components. 
\end{enumerate}
\end{remark}

\begin{ex}\label{GB04}
Let $(\Sigma,\sigma_i)$ be a stable curve of genus zero with four marked points.  
If $\Sigma$ is smooth, Gieseker bundles are  ordinary bundles, classified by their 
degree and have automorphism group $\Cs$. Assume $\Sigma$ is nodal, with two 
components $\Pn_+^1$ and $\Pn_-^1$ meeting at an ordinary double point. 
Gieseker bundles now come in two flavors:
\begin{figure}[htb]
\[\begin{xy}
(-10,0)*\xycircle(6,6){-} *{d_+};
(-10,-10)*{\mathbb{P}^1_+};
(-13,3)*{\bullet};(-13,-3)*{\bullet};
(2,0)*\xycircle(6,6){-} *{d_-};
(2,-10)*{\mathbb{P}^1_-};
(5,3)*{\bullet};(5,-3)*{\bullet};
(30,0)*\xycircle(6,6){-}*{d_+};
(30,-10)*{\mathbb{P}^1_+};
(27,3)*{\bullet};(27,-3)*{\bullet};
(42,0)*\xycircle(6,6){-}*{1};
(54,0)*\xycircle(6,6){-}*{d_-};
(54,-10)*{\mathbb{P}^1_-};
(57,3)*{\bullet};(57,-3)*{\bullet};
\end{xy}\]
\caption{The permissible modifications of $\Pn^1_+ \cup \Pn^1_-$.  
Components are labelled with the degrees.}
\end{figure}
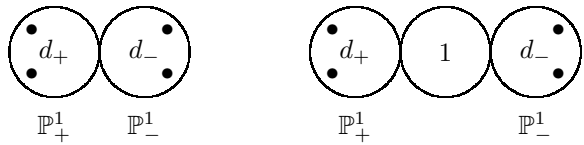
\end{ex}
\begin{enumerate}
\item  Ordinary $\Cs$-bundles on $\Sigma$, for which $m$ is the identity.  
These are classified (up to isomorphism) by their multi-degrees $\md = (d_+,d_-)$.  The 
automorphism group of any such bundle is a copy of $\Cs$, rescaling the fibers. 
\item  Gieseker bundles $(C,\mP)$ on the unique modification $m: C \to \Sigma$:  
these are also classified by their multi-degrees $\md = (d_+,1,d_-)$.  The automorphism group 
of any such Gieseker bundle is a copy of $(\Cs)^2$ ({\it cf.} Remark~\ref{automph}).
\end{enumerate}

\begin{ex}\label{GBon11}
Let $(\Sigma,\sigma_1)$ be the boundary divisor in $\Mbar_{1,1}$, representing a curve whose 
single rational component has one marked point ($\sigma_1$) and one self-node ($\sigma$).  
The degree $d$ Gieseker bundles on $(\Sigma,\sigma_1)$ come in two flavors, ordinary $\Cs$-bundles on $\Sigma$ and bundles on the modification of $\Sigma$ at $\sigma$.  
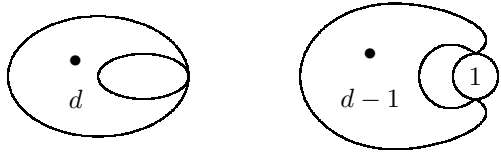
\begin{figure}[htb]
\[\begin{xy}
(40,0)*{\begin{xy}
(0,0)*\xycircle(12,8){-};
(6,0)*\xycircle(6,3){-};
(-3,-3)*{d};
(-3,2)*{\bullet};
\end{xy}};
(80,0)*{\begin{xy}
(14,0)*\xycircle(3,3){-};
(14,0)*{1};
(14,3)*{};(14,-3)*{} **\crv{(19,6) & (0,14) & (-14,0) & (0,-14)& (19,-6)};
(14,3)*{};(14,-3)*{} **\crv{(10,6) & (3,0) & (10,-6)};
(0,-3)*{d-1};
(0,3)*{\bullet};
\end{xy}};
\end{xy}\]
\caption{A pictures of $(\Sigma,\sigma_1)$ and its modification $(C,\sigma_1)$.}
\end{figure}

The normalization of $\Sigma$ is a copy of $\Pn^1$.  Line bundles on $\Pn^1$ are classified by 
their degree, so the ordinary degree $d$ bundles on $\Sigma$ are classified (up to isomorphism) 
by the gluing data at $\sigma$.  The space of such  data is a copy of $\Cs$.  
$\sigma$ is a self-node, so these gluing isomorphisms are fixed by bundle rescalings.  Hence the automorphism group of an ordinary bundle on $\Sigma$ is a copy of $\Cs$.

The modification $C$ has two rational components and two nodes, so a bundle on $C$ is specified 
by two gluing isomorphisms.  However (Remark~\ref{automph}), the automorphisms of the bundle 
over the bubble act simply transitively on the set of gluing isomorphisms, so (up to isomorphism) 
there is only {\it one} Gieseker bundle on $\widetilde{C}$.  
The automorphism group of this bundle is a copy of $\Cs$.
\end{ex}

\subsection{The Moduli Stack of Gieseker Bundles on Stable Curves}\label{universal}
Families of Gieseker bundles are classified by a moduli stack, as follows.

\begin{definition}\label{Giestack}
The {\it stack $\Mtwid_{g,I}([\pt/\Cs])$ of Gieseker $\Cs$-bundles on stable genus $g$, 
$I$-marked curves} is a fibered  category (over $\C$-schemes).  Its objects are tuplets 
$(B,C,\sigma_i,\mP)$ consisting of 
\begin{enumerate}
\item a test scheme $B$,
\item a flat projective family $\pi:C \to B$ of pre-stable, genus $g$ 
curves with marked points $\sigma_i: B \to C$, $i\in I$, and 
\item a principal $\Cs$-bundle $p:\mP \to C$ defining a family of Gieseker bundles on 
the stabilization $C\to C^{st}$.
\end{enumerate}
The morphisms in this category are commutative diagrams
\[
\xymatrix{\mP' \ar[r]^{\tilde{f}} \ar[d]^{p'} & \mP \ar[d]^p \\ 
C' \ar[r]^{f} \ar[d]^{\pi'} & C \ar[d]^\pi \\ 
B' \ar[r] \ar@/^/[u]^{\sigma'_i} & B\ar@/_/[u]_{\sigma_i}}
\]
where $\tilde{f}$ is $\Cs$-equivariant and $C' = B' \times_B C$.
\end{definition}

In section \ref{ElemGeom}, we will see that this fibered category is an Artin stack.  Below, we review the universal features of $\Mtwid_{g,I}([\pt/\Cs])$  
and discuss some examples. 

The language of stacks is designed to track automorphisms of objects in families.  For instance, the quotient stack $[\pt/\Cs]$ is the classifying stack for principal $\Cs$-bundles: any principal $\Cs$-bundle $p: \mP \to S$ on a scheme $S$ is pulled back ({\it classified}) by a unique morphism $\phi: S \to [\pt/\Cs]$ from the {\it universal} principal $\Cs$-bundle $\pt\to [\pt/\Cs]$.
 
The stack $\Mtwid_{g,I}([\pt/\Cs])$ carries several tautological families (which are 
jointly universal, subject to the Gieseker constraint). We have a family of semi-stable curves 
with marked points indexed by $I$,
\[
\pi: C \to \Mtwid_{g,I}([\pt/\Cs]), \qquad \sigma_i:  \Mtwid_{g,I}([\pt/\Cs]) \to C
\] 
and a principal $\Cs$-bundle 
\[
p: \mP \to C.
\] 
Its classifying map $\phi: C \to [\pt/\Cs]$ leads to the 
 {\it evaluation maps} $\ev_i = \phi\circ \sigma_i$
\[
\ev_i: \Mtwid_{g,I}([\pt/\Cs]) \to [\pt/\Cs] \qquad (i \in I).
\] 
Finally, there is a natural morphism
\[
F: \Mtwid_{g,I}([\pt/\Cs]) \to \Mbar_{g,I}
\]
which forgets $\mP$ and sends the curve $(C,\sigma_i)$ to its stabilization.
\subsection{Examples \& Local Model}

We treat three examples here.  The first is trivial.  The second is the simplest example in 
which the Gieseker bubbles can be seen.  The third one illustrates why our stack is typically 
non-separated and infinite type, and it displays the finite-type subspaces which we will use 
in proving our main theorem.  

\begin{ex}[genus zero, 3 marked points]
$\Mbar_{0,3}$ is a point; there is, up to equivalence, only one stable genus zero curve $\Sigma$ with 3 marked points.  The only modification of $\Sigma \simeq \Pn^1$ is the trivial modification $m=\id$.  Likewise, up to equivalence, there is one bundle of degree $D$ on $\Pn^1$, namely the principal bundle associated to $\mc{O}_{\mb{P}^1}(D)$.  The automorphism group of $\mc{O}_{\mb{P}^1}(D)$ is a copy of $\Cs$, and any $D \in \Z$ is allowed, so we conclude
\[
\Mtwid_{0,3}([\pt/\Cs]) \simeq \bigsqcup_{D\in\Z}  [\pt/\Cs].
\]
\end{ex}

\begin{remark}
The connected components of $\Mtwid_{g,I}([\pt/\Cs])$ are always classified by the total degree 
(Corollary~\ref{degreecomponents}):
\[
\Mtwid_{g,I}([\pt/\Cs] = \bigsqcup_{D\in\Z} \Mtwid^D_{g,I}([\pt/\Cs].
\]
In the next two examples we'll fix the total degree $D$.
\end{remark}
\begin{ex}[genus one, 1 marked point]\label{Ex11}
We can represent the stack of degree $D$ Gieseker bundles as a quotient 
$$
\Mtwid^D_{1,1}([\pt/\Cs]) \simeq [\mc{J}/\Cs],
$$ 
where $\mc{J}$ is the Deligne-Mumford stack which classifies triplets $((\Sigma,\sigma),(m,\mP),t)$ consisting of a genus $1$ curve $\Sigma$ with a single marked point $\sigma$, a Gieseker bundle $(m,\mP)$ of degree $D$ on a modification of $\Sigma$, and a trivialization $t: \mP_\sigma \simeq \Cs$ of the fiber of the $\mP$ at the marked point $\sigma$.  

Note that $\mc{J}$ comes equipped with a forgetful map $f:\mc{J} \to \Mbar_{1,1}$. This morphism 
is representable, the stack nature of $\mc{J}$ being captured by the base.  Thus, we obtain a 
chart $A$ for $\Mtwid^D_{1,1}([\pt/\Cs])$ by pulling $\mc{J}$ back along any chart $B \to \Mbar_{1,1}$.

Over smooth curves $(\Sigma,\sigma_1)$, the fibers of $f$ are copies of the Jacobian $\Jac(\Sigma)$ of $\Sigma$.  The Jacobian of $\Sigma$ is a copy of $\Sigma$, so over $\mc{M}_{1,1}$, the stack $\mc{J}$ is simply a copy of the universal curve $\Sigma_{1,1}$.  In fact, this is also true over the boundary divisor of $\Mbar_{1,1}$.  We saw in Example \ref{GBon11} that the space of Gieseker bundles (with a trivialization $t$ to eliminate the global rescaling automorphisms) is obtained by gluing a copy of the point $\pt$ to a copy of $\Cs$.  This gluing identifies the $\pt$ with both $0$ and $\infty$, so the resulting curve is a rational curve with a self-node.

The $\Cs$ acts on $\mc{J}$ by rescaling the trivialization $t$.  This action is easily seen to be trivial.   Thus, $\Mtwid^D_{1,1}([\pt/\Cs])$ is isomorphic to the product of the universal curve on $\Mbar_{1,1}$ and the quotient stack $[\pt/\Cs]$.

\end{ex}

\begin{ex}[genus zero, $4$ marked points]\label{Ex04}  

Recall that $\Mbar_{0,4}$ is isomorphic to $\Pn^1$.  The open locus classifying smooth curves is $\Pn^1 \setminus \{0,1,\infty\}$, and the boundary divisor $\{0,1,\infty\}$ classifies reducible nodal curves with two marked points on each smooth components.  

Let $B = \Pn^1 \setminus \{1,\infty\} = \mb{A}^1 \setminus 1$ and consider the family $(\Sigma,\sigma_i): B \hookrightarrow \Mbar_{0,4}$ of marked curves obtained by restricting the universal marked curve $\Sigma_{0,4}$ to $B$.  This family is a deformation of the curve $(\Sigma_o,\sigma_{o,i})$ which has two components meeting at a common node and each carrying two marked points.  
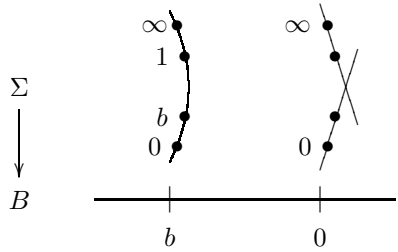
\begin{figure}[htb]
\centering
\[\begin{xy}
(0,0)*{}; (40,0)*{} **\dir{-};
(10,0)*{|};
(10,-5)*{b};
(30,0)*{|};
(30,-5)*{0};
(10,5)*{};(10,25)*{} **\crv{(15,15)};
(11,7)*{\bullet}; (8,7)*{0}; (28,7)*{0};
(11,23)*{\bullet}; (8,23)*{\infty}; (27,23)*{\infty};
(12,11)*{\bullet}; (9,11)*{b}; 
(12,19)*{\bullet}; (9,19)*{1};
(30,4)*{};(35,20)*{} **\dir{-};
(30,26)*{};(35,10)*{} **\dir{-};
(31,7)*{\bullet}; 
(31,23)*{\bullet};
(32,11)*{\bullet};
(32,19)*{\bullet};
(-10,0)*{B};
(-10,15)*{\Sigma};
{\ar@{->} (-10,12)*{};(-10,3)};
\end{xy}\]
\caption{The family $\Sigma$ over $B = \Pn^1 \setminus \{1,\infty\}$.}
\end{figure}


We will describe the stack $\Mtwid^D_{0,4}([\pt/\Cs])$ of bundles of total degree $D$ by giving a chart $A_D$ for the restriction 
$F|_B: \Mtwid \to B$ of $F:  \Mtwid_{0,4}([\pt/\Cs]) \to \Mbar_{0,4}$ to $B$.  
(Similar descriptions apply near $1$ and $\infty$ in $\Mbar_{0,4}$, and 
$\Mtwid_{0,4}([\pt/\Cs])$ is obtained by gluing these local descriptions.)  

Let $V$ be the subset $\{0,\infty\} \subset I$.  The chart $A_D$ is the algebraic space 
(scheme, in fact) which classifies isomorphism classes of tuplets
\[
(C,\sigma_i,\mP,t_0,t_\infty)
\]
where $(C,\sigma_i,\mP)$ is a family of degree $D$ Gieseker bundles on the curve $\Sigma/B$ and the $t_v \in \Gamma(B,\sigma_v^*\mP)$ are families of trivializations based at the marked points $0$ and $\infty$.  (For any closed point $b \in B$, $t_v(b)$ is a point in the fiber of $\mP$ at $\sigma_v(b)$.)  

The scheme $A_D$ carries a forgetful morphism $f: A_D \to B$, whose fibers classify trivialized bundles living on modifications of the fibers of $\Sigma \to B$.

$A_D$ has a natural stratification which classifies curves and bundles by their topological type and multi-degree.  Let $\gamma$ denote the modular graph of the smooth fibers of $\Sigma$, and let $\gamma_o$ and $\tau_o$ denote respectively the 2-vertex modular graph of $\Sigma$'s singular fiber and its unique (3-vertex) modification.

The open stratum $A_{\gamma,D}$ classifies degree $D$ bundles on smooth curves equipped with a pair of trivializations.  The global rescaling automorphisms allow us to fix one of the two trivializations, but not both, so we get an extra degree of freedom which measures the ratio of the two trivializations.  Thus, the generic fiber of the forgetful morphism $f: A_D \to B$ is isomorphic to $\Cs$.

In the special fiber of $f$, we have two kinds of strata, corresponding to $\gamma_o$ and $\tau_o$. The bundles which can appear on curves having these modular graphs were classified in Example~$\ref{GB04}$.  

For bundles on curves of type $\gamma_o$, adding the trivializations fixes the automorphisms and leaves us an extra degree of freedom, the ratio of the two trivializations, or the gluing isomorphism over the node.  The total degree $D$ can split arbitrarily between the two components, so the stratum labelled by $\gamma_o$ breaks up into $\Z$-many copies of $\Cs$.  The stratum $A_{\gamma_o, n}$ classifying twice-trivialized bundles of multi-degree $(D+n,-n)$ is a copy of $\Cs$.

Bundles on curves of type $\tau_o$ have automorphism group $(\Cs)^2$, but these automorphisms are fixed by the trivializations.  Thus, the stratum labelled by $\tau_o$ breaks up into $\Z$-many points.  We'll denote the point classifying trivialized bundles of multi-degree $(D+n-1,1,-n)$ by $A_{\tau_o,n}$.

If we smooth a Gieseker bubble, the unit of degree it carries migrates onto the component on the other side of the node.  The special fiber over $0 \in B$
$$
(A_D)_0 \simeq  \cup_{n\in\Z} \Pn^1_{n}
$$
is an endless chain of rational curves, as illustrated below. 

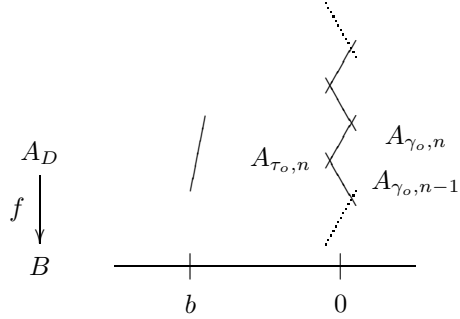
\begin{figure}[htb]\label{Achain}
\centering
\[\begin{xy}
(0,0)*{}; (40,0)*{} **\dir{-};
(10,0)*{|};
(10,-5)*{b};
(30,0)*{|};
(30,-5)*{0};
(10,10)*{};(12,20)*{} **\dir{-};
(28,3)*{};(32,10)*{} **\dir{.};
(32,8)*{};(28,15)*{} **\dir{-};
(28,13)*{};(32,20)*{} **\dir{-};
(32,18)*{};(28,25)*{} **\dir{-};
(28,23)*{};(32,30)*{} **\dir{-};
(32,28)*{};(28,35)*{} **\dir{.};
(40,17)*{A_{\gamma_o,n}};
(40,11)*{A_{\gamma_o,n-1}};
(22,14)*{A_{\tau_o,n}}; 
(-10,0)*{B};
(-10,15)*{A_D};
{\ar@{->} (-10,12)*{};(-10,3)};
(-13,7.5)*{f};
\end{xy}\]
\caption{A cartoon of $A_D$; its fibers over $B$ are drawn.}
\end{figure}

The action of the group $(\Cs)^V = \Cs_0 \times \Cs_\infty$ which rescales the trivializations preserves the forgetful morphism $f$.  The diagonal $\Cs_\Delta \subset (\Cs)^V$ acts trivially.  However, the quotient group $(\Cs)^V/\Cs_\Delta$ does not act trivially; it rescales each copy of $\Cs$ which appears in the fibers of $f$.  Its only fixed points are the points $A_{\tau_o,n}$ which classify Gieseker bundles which live on the unique, non-trivial modification of the special fiber of $\Sigma$. 
\end{ex}

\begin{notation}[Local Model]\label{localmodel}
It will be useful to have an explicit coordinate cover of the chart $A_D$ defined in the previous example. Let
\[
A_{D,n} = \Spec \C[z_n,w_n] \simeq \mb{A}^2.
\]
We make $A_{D,n}$ a coordinate chart by identifying the stratum $A_{\tau_o,n}$ with the origin, and the strata $A_{\gamma_o,n}$ and $A_{\gamma_o,n-1}$ respectively with the $z_n$ axis and the $w_n$ axis.

We cover $A_D$ by gluing $A_{D,n-1}$ to $A_{D,n}$, identifying the open sets $U_{w_n}$ and $U_{z_{n-1}}$ via the relation $w_n = 1/z_{n-1}$.  (Here $U_{w_n} = \Spec \C[z_n,w_n]_{(w_n)} \simeq \mb{A} \times \Cs$; likewise, $U_{z_{n-1}}$.)  The morphism $f: A_D \to B$ is given on $A_{D,n}$ by $b=z_nw_n$, where $b$ is the standard coordinate on $B \subset \Spec \C[b]$.  

The group $(\Cs)^2$ acts on $A_D$, with weight $(1,-1)$ on $z_n$ and weight $(-1,1)$ on $w_n$.  Thus, the diagonal acts trivially, and the fixed points are picked out by $z_n=w_n=0$.  

\end{notation}

\begin{remark}\label{Achain2}
$A_D$ has several interesting {\it finite type} subschemes of $A_D$, which we now describe.

First, consider the finite type subscheme $A^o_D$, given by
\[
A^o_D = A_{\gamma,D} \sqcup A_{\gamma_o, 0}
\]
and also for any $\delta \geq 0$, the finite-type subschemes 
\begin{align*}
Z_\delta &= A_{\tau_o,\delta+1} \sqcup A_{\gamma_o,\delta + 1} \simeq \mathbb{A}^1\\ 
W_\delta &= A_{\tau_o,-\delta} \sqcup A_{\gamma_o,-\delta - 1} \simeq \mathbb{A}^1
\end{align*}

We emphasize the following properties of these subschemes.
\begin{enumerate}
\item $A^o_D$ has no $(\Cs)^2$-fixed points, and its quotient 
\[
[A^o_D/(\Cs)^2] = B \times [\pt/\Cs_\Delta].
\]
is the product of  the classifying stack $[\pt/\Cs_\Delta]$ and a scheme $Q \mbox{ ($=B$)}$ which is proper over $B$.
\item $Z_\delta$ (resp. $W_\delta$) is a bundle of affine spaces over the fixed point locus $A_{\tau_o,\delta + 1}$ (resp. $A_{\tau_o,-\delta}$).

\item $A^o_D$ is the open locus of a stratification of $A_D$, with unions of the $Z_\delta$ and $W_\delta$ exhausting the complement $A_D \setminus A^o_D$.
\end{enumerate}
In Sections \ref{usefulatlas} and \ref{fixedpoints}, we will introduce analogues of the $A^o_D$, $Z_\delta$, and $W_\delta$ for the general case.
\end{remark}

\section{Geometry of the Stack of Gieseker Bundles}\label{ElemGeom}

We review some basic facts about the geometry of the stack of 
$\Cs$-bundles on prestable curves and its substack of Gieseker bundles. This section 
is mainly an orientation to results in the literature, so some proofs are sketched 
or omitted.

\subsection{The Stack of All Bundles on Prestable Curves}

Let $\mathfrak{C}_{g,I}\to \mathfrak{M}_{g,I}$ be the universal curve over the stack of prestable, 
genus $g$, $I$-marked curves.  

\begin{definition}
The {\it stack $\mathfrak{M}_{g,I}([\pt/\Cs])$ of principal $\Cs$-bundles on prestable $I$-marked curves of genus $g$} is the relative Hom-stack\footnote{This is not as bad as it looks.  
Note that $\mathfrak{C}_{g,I} \to \mathfrak{M}_{g,I}$ is representable.}
\[
\Hom_{\mathfrak{M}_{g,I}}(\mathfrak{C}_{g,I},[\pt/\Cs] \times \mathfrak{M}_{g,I}).
\]
 The substack of $\mathfrak{M}_{g,I}([\pt/\Cs])$ which classifies curves with modular graph $\gamma$ and bundles of multi-degree $d$ is denoted $\mathfrak{M}_{\gamma,\md}$.  
\end{definition}

\begin{proposition}
$\mathfrak{M}_{g,I}([\pt/\Cs])$ is an Artin stack, as is every substack $\mathfrak{M}_{\gamma,\md}$.
\end{proposition}
\begin{proof}[Idea of proof]
The base stack $\mathfrak{M}_{g,I}$ of prestable curves is Artin, and we can give 
Quot-scheme presentations of the stacks of bundles, locally over the base of the forgetful map 
$\mathfrak{M}_{g,I}([\pt/\Cs]) \to \mathfrak{M}_{g,I}$.
\end{proof}

\begin{proposition}
The substacks $\mathfrak{M}_{\gamma,\md}$ are of finite type and finite presentation.
\end{proposition}

\begin{proof}
We have fixed the topological type, so we may exploit the normalization of $\mathfrak{C}_{g,I}$ 
over $\mathfrak{M}_{\gamma,\md}$ to represent bundles by their lifts to the 
connected components of the normalization plus gluing data. 
\end{proof}

\begin{proposition}
The substacks $\mathfrak{M}_{\gamma,\md}$ stratify $\mathfrak{M}_{g,I}([\pt/\Cs])$:
they are locally closed and disjoint.  The whole stack is a union 
\[
\mathfrak{M}_{g,I}([\pt/\Cs]) = \bigsqcup_{\gamma,\md} \mathfrak{M}_{\gamma,\md},
\]
over all topological types.  Moreover, the closure of $\mathfrak{M}_{\gamma,\md}$ in $\mathfrak{M}_{g,I}([\pt/\Cs])$ is a disjoint union of other such strata.
\end{proposition}
\begin{proof}
This is a direct analogue of the stratification of the stack of prestable maps to projective varieties 
by labelled modular graphs, described by Behrend \& Manin in \cite{MR1412436}, Sections 1 \& 2.   
It is obtained from the standard modular graph stratification of the base stack $\mathfrak{M}_{g,I}$ 
of prestable curves by adding degree labels to track the connected components.  The modular graph stratification of prestable curves is also described in \cite{MR1412436}, but has a longer history in the literature, being implicit in Knudsen's original description \cite{K2} of contraction and clutching maps.
\end{proof}

\begin{lemma}\label{defmodgraph}
Let $(C_o,\sigma_{o,i},\mP_o)$ be a $\Cs$-bundle on a prestable curve having topological type $(\gamma_o,\md_o)$.  Suppose that we are given a deformation $(C,\sigma_i,\mP)$ of $(C_o,\sigma_{o,i},\mP_o)$ over the Spec of a complete discrete valuation ring.  The topological type $(\gamma,\md)$ of the generic fiber can be any degree-labelled modular graph obtained from $(\gamma_o,\md_o)$ by finite combinations of the following elementary operations:
\begin{enumerate}
\item Resolve a self node: Delete a self-edge attached to a vertex $v$, increasing the genus $g_v$ by 
$1$, leave the multi-degree unchanged.
\item Resolve a splitting node:  join a pair of adjacent vertices $v_1$ and $v_2$ into a single 
vertex $v$, having genus $g_v = g_{v_1} + g_{v_2}$ and degree $d_v = d_{v_1} + d_{v_2}$.  
Delete one edge joining $v_1$ and $v_2$, and convert the others to self-edges.
\end{enumerate}
Moreover all such modular graphs occur in some deformation.

\end{lemma}
\begin{proof}
First note that any deformation of $C_o$ over a complete DVR will only smooth nodes; new nodes 
are not created.  This limits the topological types $\gamma$.  To determine what 
multi-degrees are allowed, observe that we may normalize $C$ at the nodes of $C_o$ which are 
not smoothed by deformation, and consider each connected component separately.  Flat 
deformations preserve the total degree on connected curves, which fixes the 
degrees which can appear.
\end{proof}

\begin{proposition}  
$\mathfrak{M}_{g,I}([\pt/\Cs])$ is locally of finite type and locally finitely presented.
\end{proposition}

\begin{proof}
The previous lemma implies that we only reach finitely many strata by deformation.  
Each such stratum is of finite type and finite presentation.
\end{proof}

\begin{proposition}
$\mathfrak{M}_{g,I}([\pt/\Cs])$ is unobstructed.
\end{proposition}
\begin{proof}
Formal deformations of the curve and bundle pair are controlled by the Atiyah complex,  
a short exact sequence combining the adjoint bundle $\mathrm{ad}\,\mc{P}$ and the 
tangent bundle of $C$ (based at the marked points)
\[
0\to \mathrm{ad}\,\mc{P} \to \mc{D} \to \mc{T}_C(-\sum\sigma_i) \to 0.
\] 
A key detail is that $\mc{T}_C$ is a line bundle (not a complex), even when 
$C$ is nodal. Deformations are tangent to $H^1(C;\mc{D})$, while obstructions live in $H^2$, 
which vanishes  because $C$ is one-dimensional.
\end{proof}

\begin{proposition} 
The dimension of $\mathfrak{M}_{g,I}([\pt/\Cs])$ is $(g-1) + 3(g-1) + |I|$.
\end{proposition}

\begin{proof}
$\mathfrak{M}_{g,I}([\pt/\Cs])$ is unobstructed, so its dimension is its virtual dimension, 
which is the negative of the Euler characteristic of $\mc{D}$.
\end{proof}

\subsection{$\Mtwid_{g,I}([\pt/\Cs])$ as a Substack of $\mathfrak{M}_{g,I}([\pt/\Cs])$}

The Gieseker condition being  topological, several facts about the stack 
of Gieseker bundles follow from corresponding facts about $\mathfrak{M}_{g,I}([\pt/\Cs])$.

\begin{proposition}
If, in Lemma \ref{defmodgraph}, $(\gamma_o,\md_o)$ satisfies the Gieseker 
conditions, then so does $(\gamma,\md)$.

\end{proposition}
\begin{proof}
The unstable vertices in $\gamma$, where the condition must be checked, are identified 
with a subset of those of $\gamma_o$.  
\end{proof}

\begin{corollary} 
$\Mtwid_{g,I}([\pt/\Cs])$ is an open substack of $\mathfrak{M}_{g,I}([\pt/\Cs])$, and 
a union of strata. In particular, 
\begin{enumerate}
\item $\Mtwid_{g,I}([\pt/\Cs])$ is an Artin stack. 
\item $\Mtwid_{g,I}([\pt/\Cs])$ is unobstructed. 
\item $\dim \Mtwid_{g,I}([\pt/\Cs])=(g-1) + 3(g-1)+|I|$. \qed
\end{enumerate}
\end{corollary}

\begin{notation}\label{degeneracystrata}
If a topological type $(\gamma,\md)$ satisfies the Gieseker condition, the 
corresponding substack of $\Mtwid_{g,I}([\pt/\Cs])$ will be denoted by 
$\mc{M}_{\gamma,\md}$.
\end{notation}

\begin{proposition}
$\Mtwid_{g,I}([\pt/\Cs])$ inherits a topological type stratification from $\mathfrak{M}_{g,I}([\pt/\Cs])$.   The substacks $\mc{M}_{\gamma,\md}$ are locally-closed and disjoint, and the whole moduli stack is their union:
\[
\Mtwid_{g,I}([\pt/\Cs]) = \bigsqcup_{\gamma,\md} \mc{M}_{\gamma,\md}.
\]
Finally, the closure of any given $\mc{M}_{\gamma,\md}$ in $\Mtwid_{g,I}([\pt/\Cs])$ is obtained as a union 
$$
\op{cl}\big(\mc{M}_{\gamma,\md}\big) = \bigsqcup_{\gamma',\md'} \mc{M}_{\gamma',\md'},
$$
where the  union is over all multi-degree-labelled modular graphs $(\gamma,\md')$ obtained from $(\gamma,\md)$ by finite combinations of the following elementary operations:
\begin{enumerate}
\item Self node: Lower the genus of a vertex by 1 and add a self-edge.
\item Splitting node:  Split a vertex $v$ into two vertices $v_1$ and $v_2$, connected by an edge, 
with $g_{v_1} + g_{v_2} = g_v$ and $d_{v_1}+ d_{v_2} = d_v$.
\item Gieseker bubbling:  Replace an edge connecting a stable vertex $v$ to a stable vertex $v'$ 
with two edges connected to a new common  vertex having $g = 0$ and $d = 1$, while subtracting $1$ 
from either $d_v$ or $d_{v'}$. (Note that $v$ may equal $v'$.) \qed
\end{enumerate} 
\end{proposition}

\begin{corollary}\label{degreecomponents}
The connected components of $\Mtwid_{g,I}([\pt/\Cs])$ are labelled by total degree $D$.  
\[
\Mtwid_{g,I}([\pt/\Cs]) = \bigsqcup_{D \in \Z} \Mtwid_{g,I}^D([\pt/\Cs]).
\]
\end{corollary}
\begin{proof}
Any Gieseker bundle may be deformed to a bundle on a smooth curve, so all Gieseker bundles with the same total degree lie in the same connected component.  Conversely, no deformation can change the total degree.
\end{proof}

\subsection{Limits of Bundles}\label{sec2}

The Gieseker stack $\Mtwid_{g,I}([\pt/\Cs])$ has infinitely many connected components, 
and even its connected components generally have infinite type: a modular graph $\gamma$ 
with at least two vertices carries countably many multi-degrees $\md: V_{\gamma} \to \Z$ 
for which $\sum d_v = D$. In addition, $\Mtwid_{g,I}([\pt/\Cs])$ is not separated, because 
of the continuous automorphism groups of line bundles. Even if we fix the fiber of the 
bundle at some marked points, Gieseker bubbles introduce additional automorphisms, which 
keep our stack typically non-separated. Nonetheless, $\Mtwid_{g,I}([\pt/\Cs])$ 
satisfies the \emph{valuative criterion for completeness}.

Let then $R$ be a complete discrete valuation ring with fraction field $K$ and   
denote by $D$ the disc $\Spec(R)$ and  by $D^\times$ the punctured disc 
$\Spec(K)$. Let $C^\times\to D^\times$ be a family of marked, pre-stable curves 
carrying a $\Cs$-bundle $\mP^\times\to C^\times$. We omit the marked points from the 
notation, and will at times impose additional restrictions on $C^\times, \mP^\times$.

\begin{proposition}\label{completeness}
Any family $(C^\times,\mP^\times)$ of Gieseker bundles can 
be extended to a Gieseker family $(C,\mP)$ over $D$ (possibly after \'{e}tale base 
change on $D^\times$).
\end{proposition}

\begin{proof} 
By completeness of $\Mbar_{g,I}$, there is  (after \'etale refinement of 
$D^\times$) a unique stable curve $\Sigma$ extending the stabilization of 
$C^\times$.  The claim now follows from the rank $1$ case of Theorem 2 of 
Nagaraj-Seshadri \cite {MR1687729}: the coarse moduli space of rank $n$ Gieseker 
bundles on the family $\Sigma/D$ is projective, hence complete.  (Nagaraj \& Seshadri 
state the result in the special case that the generic fiber of $C$ is smooth and 
the special fiber irreducible with a single self-node; however, their argument for 
the existence of limits for the bundles, in Section~4 and the appendix of 
\cite {MR1687729}, is local near the nodes, so the more general case  follows.)
\end{proof} 

\begin{remark}
Uniqueness fails in this Proposition, but with discrete ambiguity: 
a choice of modification $C$ of $\Sigma$ (with fixed $C^\times$), and the 
extension of $\mP^\times$ thereon. Now, two extensions $\mP, \mP'$ of $\mP^{\times}$ 
to a given $C$ differ by a 
{\it twister},  a line bundle represented by a Cartier divisor $\sum_{v \in 
V_{\gamma}} m_v C_{o,v}$ on $C$. In the special case when the total space of $C$ is 
regular, Caporaso \cite{MR2382140} shows that twisters are determined by their 
multi-degree on $C_o$, and identifies all possibilities.

If $(C,\mP)$ is a Gieseker bundle, contracting the bubbles leads to 
$A_1$-singularities in $\Sigma_o$, so the location of the bubbles is pinned by 
$\Sigma$. In that case, all information about $\mP$ is contained in the multi-degree 
of $\mP$ over $C_o$.
\end{remark}

\section{Local Presentation of the Gieseker stack}\label{usefulatlas}

We describe a local quotient presentation $A/\mc{G}$ for our stack of Gieseker bundles.  
Then we refine the topological type stratification to one which tracks the nodes 
being  smoothed in a deformation.  

\begin{notation}
We fix the following notation for this section and the next.
\begin{enumerate}
\item $(\Sigma_o,\sigma_{o,i})$ is a stable marked curve of type $(g,I)$
\item $V = V_{\gamma_o}$ denotes the set of vertices of its modular graph $\gamma_o$
\item $d:V\to\Z$ is a general multi-degree, giving a topological type $(\gamma_o,d)$
\item $\mc{G}$ is the group $(\Cs)^V$, $\Cs_\Delta\subset\mc{G}$ the diagonal subgroup
\item $B$ will be an affine \'{e}tale neighborhood of $(\Sigma_o,\sigma_{o,i})$ 
in $\Mbar_{g,I}$. It carries a locally universal deformation $(\Sigma,\sigma_i)$ of 
$(\Sigma_o,\sigma_{o,i})$. 
\item $\sigma_v: B \to \Sigma$ ($v \in V$) is an additional set of smooth 
special points over $B$, each meeting the respective component $v$ of $\Sigma_o$.  
We also assume that every stable component of every fiber of $\Sigma$ carries a $\sigma_v$.  
\item $\sigma_+$ is a particular chosen $\sigma_v$.
\end{enumerate}
\end{notation}

\begin{remark}
We may need to shrink $B$ repeatedly in later discussion, but 
in any case, the stack $\Mbar_{g,I}$ can be covered by finitely many desirable $B$s.
\end{remark}

\subsection{The space $A$}\label{sec1}
Denote by  $\Mtwid|_B$ the fiber of $F$ over $B\to \Mbar_{g,I}$ under the forget-and-stabilize 
map
\[ 
F: \Mtwid_{g,I}([\pt/\Cs]) \to \Mbar_{g,I}.
\]
We will present $\Mtwid|_B$ as a quotient stack by trivializing the bundles at special points. 
Recall that for a prestable curve $C$ with a special point $\sigma: B \to C$, a {\it 
trivialization  at $\sigma$} of a principal $\Cs$-bundle $\mP\to C$ is an 
isomorphism $t: \sigma^*\mP \to \Cs$ with the trivial bundle over $B$. 
(We may need to refine $B$ for $t$ to exist.) 

\begin{definition}\label{prop_atlas}
The {\it local chart} $A$ for $\Mtwid|_B$ is the stack of Gieseker bundles over the 
curve $\Sigma\to B$, equipped with a trivialization $t_v$ at each $\sigma_v$;  
isomorphisms are required to be compatible with the trivializations.  

Denote by $A_D\subset A$ the connected component of bundles of total degree $D$. 
\end{definition}

\begin{remark}
This stack is a category fibered over $B$, whose formal definition follows the 
template of Definition~\ref{Giestack}. The reader is entrusted to write 
out all ingredients of objects and morphisms, minding that morphisms must 
preserve all the structure, and the test scheme $X$ must now live 
over (our fixed) $B$. 
\end{remark}

\begin{proposition} \label{stackisspace}
The stack $A$ is represented by an algebraic space.
\end{proposition}
\begin{proof}
It is enough to check that the geometric points of $A$ have no automorphisms. 
Fix therefore a $\Sigma$ and $\mP\to C$. By Remark~\ref{automph}, $\Aut(\mP)$ is 
computed by deleting the Gieseker bubbles from $C$. This, however, leaves the 
stable components, each of which carries at least one trivialization point for $\mP$. 
\end{proof}

The group $\mc{G} =(\Cs)^V$ acts on $A$ by scaling the trivializations, and  
displays $\Mtwid|_B$ as a quotient stack
\[
\Mtwid|_B = A/\mc{G}.
\]

\begin{corollary}
$A$ and $\Mtwid_{g,I}([\pt/\Cs])$ are smooth.
\end{corollary}

\begin{proof}
$\Mtwid_{g,I}([\pt/\Cs])$ is unobstructed, hence formally smooth, and is 
locally of finite presentation. Thus, $A$ is formally smooth and locally 
of finite presentation, therefore smooth. Finally, the quotient of a smooth 
algebraic space by a smooth group action is a smooth stack.
\end{proof}

\subsection{The Stable Subspace $A^o$}\label{stable}

For each total degree $D$, we will identify an open subspace $A^o_D \subset A_D$ for which the quotient stack $[A^o_D / \mathcal{G}]$ is the product of the stack $[\pt/\Cs_\Delta]$ and a smooth proper quotient space $Q/B$.  The union of the $A_D$ is $A^o$.

Twisting line bundles with our chosen preferred point $\sigma_+$ equivariantly identifies the various spaces $A_D$.  We'll define $A^o_G$ first, where the total degree is the genus $G=g(\gamma_o)$, and then extend the definitions using these isomorphisms.

\begin{remark}  $A^o_G$ classifies multiply trivialized bundles in the substack $\Mtwid^o\subset\Mtwid_{g,I}([\pt/\Cs]$ of Caporaso's {\it stably balanced bundles} \cite{MR2382140}, the case when degree = genus being one of the favorable ones for her construction of the universal Picard stack.  (The bounds in Definition~4.6 of {\it loc.~cit.} are not integral and therefore strict.)  We could cite her results and skip the rest of this section, but we are including our brief treatment to keep the discussion self-contained.  
\end{remark}

\begin{definition}\label{bounddef}
A Gieseker bundle $(\mP,C)$ of total degree $G$ {\it meets the genus bounds} if 
its restriction to any subcurve $S\subset C$ has degree no less than the genus $g(S)$.   (Likewise, a Gieseker bundle $(\mP,C)$ of total degree $D$ meets the genus bounds if the $(G-D)\sigma_+$ twist of $\mP$ does.)

The {\it stable subspace} $A^o_G \subset A_G$ comprises the bundles meeting the genus bounds.  $A^o_D$ is the appropriate twist of $A^o_G$.
\end{definition}

\begin{proposition}
$A^o_G\subset A_G$ is open and of finite type.
\end{proposition}
\begin{proof}
The genus bounds are conditions on the topological type. The elementary operations 
of Lemma~\ref{defmodgraph} preserve the genus bounds, proving openness. 
Next, an upper bound on the degree on each $S$ follows from the lower bound on 
its complement, so $A^o_G$ is a union of finitely many topological type strata.  
\end{proof}

\begin{remark}
Equivalent bounds are enforced by the collection of inequalities 
\begin{equation}\label{genusbound2}
\deg(\mP|_S) > g(S)-h^0(S), \quad\forall S\neq\emptyset
\end{equation}
which need testing only against connected $S$.
Additivity makes the right-hand side, at times, more convenient than $g$ alone.
\end{remark}

\begin{ex}\label{bubbletree} Here are some illustrations of the genus bounds:
\begin{enumerate}
\item Attaching a new Gieseker bubble to $C$, at two arbitrarily specified marked 
points, preserves the genus bounds; so does erasing any existing bubble.
\item A {\it tree} meets the genus bounds iff $d_v = g_v$ on each 
component $v$: since $\sum d_v =\sum g_v$, strict inequality somewhere would break 
the genus bound elsewhere. In particular, all components are stable.  
\item Decorating a such tree arbitrarily with Gieseker bubbles as in (1) gives 
more examples, which we call {\it  Christmas trees}. (They are not trees.)
\end{enumerate}
\end{ex}

\begin{proposition}
The Gieseker topological types appearing in $A^o_G$ are precisely the deformations of 
Christmas trees (\ref{bubbletree}.3) of total degree $G$. 
\end{proposition} 
\begin{proof}We just need to find a Christmas tree degeneration for a connected 
Gieseker type graph meeting the genus bound. 

Deleting a Gieseker bubble cannot disconnect the curve, else one of the resulting 
components would break the genus bound (we have lost total degree but no genus); 
so we start by deleting all existing Gieseker bubbles. If we don't have a tree yet, 
it suffices to produce a single new Gieseker bubble by degeneration. Call a full 
subgraph $\eta\subset\tau$ 
{\it strict} if $\deg(\eta) = g(\eta)$; the full graph $\tau$ is strict. 
Observe that the intersection of two strict subgraphs is also strict: any excess $d>g$ 
on the intersection would lead to a deficit on the union of the graphs, so certainly on 
the full subgraph they span. 

Every vertex thus lies in a minimal strict subgraph. If every vertex is 
strict and there are no bubbles, then the graph is already a tree (there is no degree 
to spare for the extra genus coming from a circuit). Else, choose a vertex $v$ with 
$d_v>g_v$, and bubble off one unit of degree along an edge $e$ within its strict
minimal subgraph. Since every strict subgraph containing $v$ also contains $e$, 
it will not break the genus bounds, while the other subgraphs containing $v$ 
can bear the loss of one unit of degree.
\end{proof}

\subsection{The Proper Quotient $A^o_D/\mc{G}$}
Even if we ignore the trivializations, the Gieseker bundles in $A^o_D$ only allow 
global rescalings as bundle automorphisms, because their Gieseker bubbles may not 
disconnect the curve. Factoring $\mc{G}\simeq\Cs_\Delta \times (\Cs)^{V \setminus \{v_+\}}$ 
leads to an equivalence
\[
[A^o_D/\mc{G}] \simeq Q \times [\pt/\Cs_\Delta],
\]
where $Q$ is an algebraic space which classifies the Gieseker bundles of 
degree $g(\gamma_o)$, trivialized at $\sigma_+$, and which satisfy the genus bounds.

\begin{theorem}\label{lbcoherence}
The  moduli space $Q$ is proper over $B$.
\end{theorem}

\begin{proof}
We will use the valuative criterion, after some simplification. As usual, 
$D, D^\times$ denote the formal disk and its punctured version, $\Sigma\to D$ 
a stable curve and $\mP^\times\to \Sigma^\times$ the family of $\Cs$-bundles 
over $D^\times$ whose Gieseker completions $(C,\mP)$ over $\Sigma$ we seek. 
We have assumed the stability of $C^\times = \Sigma^\times$, since  Gieseker 
bubbles may be attached and removed at will.

We will make additional restrictions, which we can afford thanks to the normal 
crossing nature of the deformation type stratification of $A$ and the compatible 
toric action of $\mc{G}$.

\medskip
{\it Existence:} We ask that $C$ develop a single new node in the special fiber. 
Since $A/\mc{G}$ is complete over $B$, $A^o$ only fails the existence test if we find 
a $D^\times$-family for which each Gieseker completion specializes in some stratum of 
$A\setminus A^o$. If so, there is such a missing stratum of highest dimension, and we 
can detect its absence by approaching it transversally in $A$, from within a next-higher 
stratum from $A^o$. This creates exactly one ordinary double point in the special fiber. 

Call $S^o\subset A^o_D$ the stratum of $(\Sigma^\times, \mP^\times)$, with topological 
type $(\gamma,d)$ meeting the genus bounds. We seek a codimension $1$ boundary 
stratum $\partial S^o \subset A_D$, containing a Gieseker limit of $\mP^\times$, which 
also meets the genus bounds. The quest is easy when $\gamma$ is a tree: 
the new node splits $\gamma$ in two, and we adjust the two incorrect degrees ---  
adjacent to the new node --- by twisting with a suitable multiple (excess over the genus) 
of one entire side of the node. 

For general $\gamma$, if luck had it that the specialization $C_o$ created one of the 
Gieseker bubbles appearing in a Christmas tree degeneration of $(\gamma, d)$, then $C_o$ 
would also meet the genus bound, and we would be done. If not, let us travel parallel to the 
stratum of any specialization $\mP_o$, so as to reach a Christmas tree degeneration 
of $(\gamma, d)$. (This journey happens in the fibers of $F:A\to B$.) At a destination
$(T^\times, \mP^\times)$, we may remove all Gieseker bubbles before classifying 
specializations. The remaining tree degenerates into a (unique) limit stratum $S_\infty
\subset A^o$ meeting the genus bounds. 
Deform back from $S_{\infty}$ by unbubbling the Christmas tree (but leaving the new node 
intact); this keeps the genus bounds and identifies our stratum $\partial S^o$ 
parallel to the original choice.
 
\medskip
{\it Uniqueness:} Here, we insist on regularity of the total space $C$. We can do 
this, because uniqueness holds for limits in the bulk stratum of smooth curves: 
multiple limits would appear when approaching boundary strata which are distinct 
in $A^o$, but whose neighborhoods are identified in the bulk. We could then detect 
this by moving into the largest such strata from  the bulk of $A^o$, 
transversally to all divisors cutting out the problem strata. (The normal structure 
is a collection of $\C^m$s identified together everywhere except at the origin, 
and we would travel diagonally.)  This implies precisely the regularity of the total 
space of $C$ (although not of $\Sigma$, which depends on the projected path in $B$). 
The $\mc{G}$-action relates all the completed families of $\mP$ and preserves transversality, 
so regularity will apply in all Gieseker completions. 

Now, the the Gieseker bubbles in $C$ are located over the ($A_1$)
singularities of the total space of $\Sigma$. Any two Gieseker extensions $\mP,\mP'$ 
thus live on the same family $C$, the minimal resolution. Thus, $\mP'= \mP(\sum a_k E_k)$, 
the twist by a sum of increasing effective divisors $E_1 \subset E_2\subset\dots$, 
with coefficients which we may take to be positive, because $C_o\sim 0$. If $\mP$ 
obeys the genus bounds, 
we claim that $\mP'$ breaks them on the support of $E_1$. Indeed, the first twist will 
lower the degree of $\mP$ at least by the number $s$ of nodes splitting $\sup(E_1)$ 
from its complement $K$ in the special fiber of $C$; but this already breaks the 
genus bound, because 
\[
\deg(\mP|\sup(E_1)) \le G -\deg(K) \le G- g(K) \le g(\sup(E_1)) +s-1.
\]
Further twisting lowers the degree even more.
\end{proof}

\section{Stratification of $A$}\label{fixedpoints}

The stable subspace $A^o_D\subset A_D$ is the open component in a stratification we will use to prove our main theorem. To define it and study its properties, 
we first introduce a {\it deformation type} stratification of $A_D$, refining the stratification by topological type. We then define spaces $Z_\delta(\pi)$, $W_\delta(\pi)$ which retract to their fixed point sets under distinguished $\Cs$-actions.

\subsection{Stratifying $A$ by Deformation Type}\label{deftypestrat}

Let us describe the deformation type stratification for $\Mbar_{g,I}$ first. 
Each stratum $\mathcal{M}_\gamma \subset \Mbar_{g,I}$ is an intersection of normally crossing branches of divisors. 

In a sufficiently small \'{e}tale chart, the branches become distinct connected 
components, each of them representing a persistent node in the deformation. 
For instance, the one-dimensional stratum in $\Mbar_{2,0}$ representing curves with 
two nodes is the self-intersection of the boundary divisor; but the double 
cover defined by labeling the nodes 
is locally an intersection of two separate divisors. 

\begin{definition}
A {\it deformation map} $c: \gamma \to \gamma'$ of modular graphs is a continuous map 
$|\gamma|\to |\gamma'|$ which sends vertices to vertices and tails to tails, while 
possibly contracting edges to vertices.   The map $c$ induces a genus labelling on 
$\gamma'$:  $g_{\gamma'}(v) = \sum_{v' \in c^{-1}(V_{\gamma'})} g(v') + \dim H^1(|c^{-1}(v)|)$.)
\end{definition}

The strata near $\mathcal{M}_\gamma$ are in one-to-one correspondence with 
deformation maps whose domain is $\gamma$.  More precisely, Lemma~\ref
{defmodgraph} gives:
\begin{proposition}
After \'{e}tale refinement, the modular graph stratification 
$B = \bigsqcup_{\gamma} B_{\gamma}$ inherited from $\Mbar_{g,I}$ refines to a 
stratification $B = \bigsqcup_{c:\gamma_o \to \gamma} B_{c}$, 
labelled by deformation maps $c: \gamma_o \to \gamma$ of the modular graph $\gamma_o$. \qed
\end{proposition}

We lift this stratification from $B$ to $A_D$, account for Gieseker bubbling, and track degrees.

\begin{proposition}\label{refineddeg}
The topological-type stratification of $A_D$ by degree-labelled modular graphs refines to a stratification $A_D = \bigsqcup_{\dt} A_{\dt}$ with labels $\dt:=(c,\tau,d)$ consisting of a deformation map $c:\gamma_o \to \gamma$, the graph $\tau$ of a modification of a curve with modular graph $\gamma$, and a multi-degree $\md: V_\tau \to \Z$. \qed
\end{proposition}

We  call $c: \gamma_o \to \gamma$ the {\it deformation type} (with respect to 
$\Sigma_o$) of any curve $\Sigma$ parametrized by $B_{c}$. 
Likewise, the {\it deformation type of a Gieseker bundle} $\mP$ 
of multi-degree $d$ on such a curve\footnote{The terminology is slightly abusive, 
because we do not track the deformation of $\mP$ from the modification $C_o\to\Sigma_o$.} 
is the triplet $\dt = (c,\tau,d)$.

\subsection{The Strata $Z$ and $W$}\label{sec3}

Now we define the promised subspaces $Z_\delta(\pi)$ and $W_\delta(\pi)$.

\begin{notation}\label{not4}
Let $\Pi(V)$ be the power set of $V$, and fix a vertex $v_+ \in V$. 
\begin{enumerate}
\item A $\pi\in \Pi(V)$ and its complement induce a pair $(\pi_+,\pi_-)$ of 
full subgraphs of $\gamma_o$, labelled so that $v_+\in\pi_+$.
\item The {\it $\pi$-splitting edges} in $\gamma_o$ are those joining 
$\pi_+$ with $\pi_-$.  
\item $\dt=(c,\tau,d)$ denotes a general deformation type of $(\gamma_0,d)$.
We say it is {\it compatible} with $\pi$ if $c(\pi_+)\cap c(\pi_-)$ are disjoint.
\item The $\dt$-deformation $\dt(\pi)$ is the  subgraph of  
$\tau$ induced by $c(\pi)\subset\gamma$ plus its {\it internal} Gieseker bubbles 
(those attaching only to $c(\pi)$).
\item The {\it genus} $g(\pi_\pm)$ is the arithmetic genus of $\pi_\pm$. 
\item $k(\pi_\pm)$ is the number of connected components of $\pi_\pm$.
\end{enumerate}
\end{notation}

\begin{remark}
If $\pi$ and $\dt$ are compatible, then $g(\pi_\pm)$ and $k(\pi_\pm)$ are 
not changed by $\dt$.
\end{remark}

\begin{definition}
Given a deformation type $\dt$ with total degree $D=G$, define the {\it defect} of a full subgraph $\eta\subset \tau$ as $\dct(\eta):= g(\eta)-k(\eta)-\deg(\eta)$, {\it if this number is non-negative}.  For general $D$, we use the chosen isomorphism $A_D \simeq A_G$ to relabel the deformation types, and then take the defect.
\end{definition}
\noindent
Thus, $\eta$ has a defect iff it breaks the genus bounds \eqref{genusbound2}. Zero is a defect. If $\tau$ meets the genus bounds, the only defective subgraph is $\emptyset$, with defect zero. 

\begin{proposition}
For the span  $\eta_1*\eta_2$ of two full subgraphs $\eta_{1,2}\subset\tau$,
\[
\dct(\eta_1*\eta_2) \ge \dct(\eta_1) + \dct(\eta_2) 
	- \dct(\eta_1\cap\eta_2),
\]
with equality only if all new edges in $\eta_1*\eta_2$ carry Gieseker bubbles.
\end{proposition}
\begin{proof}
This follows from the additivity of $d, g, k$, accounting for  new edges.
\end{proof}
By considering all subgraphs achieving the maximal value of the defect and repeatedly 
applying the proposition, we find 
\begin{corollary}
Every deformation type $\dt=(c,\tau,d)$ contains a largest full subgraph 
$\mu(\dt)\subset\tau $ of maximal defect.\qed
\end{corollary}

This graph contains all of its internal Gieseker bubbles, so is determined 
by the stable subgraph in $\gamma$, or its $c$-preimage in $\gamma_o$. 
Comparing the various $\mu(\dt)$ within $\gamma_o$, we note a lexicographic 
ordering on pairs $(\dct(\mu(\dt)), \mu(\dt))$, which is compatible with deformations: 
if $\dt$ deforms to $\dt'$, then $\dt\ge\dt'$, and equality is 
only preserved by deformations which do not smooth any nodes joining $\mu$ with its 
complement 
in $\tau$.

\begin{definition}\label{defzw}
Given $\delta\ge 0$ and $\pi\in\Pi(V)$, let $T(\pi,\delta)$ 
be the collection of $\pi$-compatible deformation types $\dt$ with maximal defect $\delta$.  Define the substacks of $A_D$ in terms of $\delta$ and the 
maximal defective subgraph $\mu$, 
\begin{align*}
Z_{\delta}(\pi)  := \bigsqcup A_\dt\left|\: \dt\in T(\pi,\delta)\text{ and } 
\mu(\dt)= \dt(\pi_-) 
\right. \\
W_{\delta}(\pi)  := \bigsqcup A_\dt\left|\: \dt\in T(\pi,\delta)\text{ and } \mu(\dt) = 
\dt(\pi_+) \right.
\end{align*}

\begin{remark}
We suppress the degree $D$ from the notation, since we'll only ever consider one $A_D$ at a time.  Note that the defect labels $\delta$ are $D$-dependent.
\end{remark}
\end{definition}

\begin{remark}\label{easy}
Here are some easy consequences of the definition.
\begin{enumerate}
\item The $Z_\delta(\pi), W_\delta(\pi)$ are pairwise disjoint and exhaust $A_D\setminus A^o_D$. 
\item Each $Z, W$ is a finite union of deformation type strata $\dt$; in particular, it is of finite type. (We get bounds on the multi-degree.)
\item The boundary of $Z_\delta(\pi)$ meets only those $Z$ and $W$ with greater defect or larger maximally defective subgraph. (Likewise for $W$.) In particular, $Z_\delta, W_\delta$ are locally closed, and together with $A^o_D$ they assemble to a stratification of $A_D$.  
\end{enumerate}
\end{remark}

\subsection{Stabilizers and Fixed Points}

A chosen $\pi \in \Pi(V)$ splits the group $\mc{G} = (\Cs)^{V}$ into $\mc{G}_+\times \mc{G}_-$. 
Call 
\[
\mc{G}(\pi):=\Delta(\mc{G}_+)\times\Delta(\mc{G}_-) \subset \mc{G}
\]
the product of the diagonal subgroups, which we also denote by $\Cs_+\times \Cs_-$. 
The proof of the following variant of Prop.~\ref{stackisspace} is left as an exercise.

\begin{proposition}\label{splitedges}
The fixed-point set $F(\pi)$ of $\mc{G}(\pi)$ on $A_D$ is the union of strata $A_\dt$  
with $\dt$ compatible with $\pi$, and for which, in addition, every splitting edge  
carries a Gieseker bubble. \qed
\end{proposition}

\begin{remark}\label{fixdegrees}
$F(\pi)$ decomposes by degree into {\it closed} subspaces; the bi-degrees 
$n_\pm:= \deg\dt(\pi_\pm)$ are constant within a connected component, because 
deformations within $F(\pi)$ may not smooth any of the Gieseker bubbles splitting $\pi_\pm$. Also, $D= s(\pi)+ n_+ + n_-$, so that $n_+$ and $n_-$ determine each other.   (Here, $s(\pi)$ is the number of $\pi$-splitting edges.)
\end{remark}

\begin{ex}
In Example \ref{Ex04}, $\pi = \{ \{v_+\}, \{v_-\}\}$ is the only non-trivial partition for $\gamma_o$.  The fixed point stratum $F_n(\pi) \in A_D$ labelled by 
$\un = (n_+,n_-)$ is the point $A_{\tau_o,n_+ + 1 - D}$ which classifies Gieseker 
bundles with-trivializations of multi-degree $(n_+,1,n_-)$. 
\end{ex}

\begin{remark}\label{generalstab}
One can describe the stabilizers of all deformation strata as follows. A deformation type $(c,\tau,d)$ induces a partition $P$ of $V$, from the connected 
components left in $\tau$ after removing all Gieseker vertices. The stabilizer on the stratum $A_{c,\tau,d}$ is the multi-diagonal subgroup $\mc{G}(P)$ for this partition.
\end{remark}

\subsection{$Z$ and $W$ as Bundles over their Fixed Point Strata}

Let $F_{\delta,z}(\pi) := F(\pi)\cap Z_\delta(\pi)$ and $F_{\delta,w}(\pi) 
:= F(\pi)\cap W_\delta(\pi)$. These are open substacks of $F(\pi)$, as per 
the discussion preceding Definition~\ref{defzw}.   

Smoothing the $\pi$-splitting Gieseker bubbles in $F_{\delta,z}$ into the 
$\pi_+$-side (likewise, in $F_{\delta,w}(\pi)$ on the $\pi_-$-side) does not 
change the maximal property of $\pi_-$, respectively $\pi_+$. More precisely, 

\begin{proposition}{\hspace{1pt}}
\par
\begin{enumerate}
\item $Z_\delta(\pi)$ classifies bundles which arise from $F_{\delta,z}(\pi)$ by smoothing 
away nodes attaching components of $c(\pi_+)$ to splitting Gieseker bubbles. 
\item $W_\delta(\pi) $ classifies bundles which arise from $F_{\delta,w}(\pi)$ by smoothing 
away nodes attaching components of $c(\pi_-)$ to splitting Gieseker bubbles.
\item $\Cs_+$ acts with positive weights on the conormal bundle of $F_{\delta,z}(\pi)$ 
in $Z_\delta(\pi)$ and on the normal bundle to $Z_\delta(\pi)$ in $A_D$. 
\item $\Cs_+$ acts with negative weights on the conormal bundle of $F_{\delta,w}(\pi)$ 
in $W_\delta(\pi)$ and on the normal bundle to $W_\delta(\pi)$ in $A_D$.
\end{enumerate}
\end{proposition}
\noindent
Recall that the diagonal in $\Cs_+\times \Cs_-$ acts trivially on $A$, so the opposite 
holds for the $\Cs_-$-weights. Fix $\pi$ for now and drop it from the notation where possible. 

\begin{proof}
The first two statements witness that the bundles in $Z,W$ develop Gieseker bubbles at 
the $\pi$-splitting nodes when we scale the $v_+$-trivialization. 

More precisely, scaling the trivialization of $\mP$ at $v_+$ to $0$ gives a map 
$\eta_{0,l}: Z_\delta \to F_{\delta,z}$, and scaling to $\infty$ gives 
$\eta_{\infty,l}: W_\delta \to  F_{\delta,w}$. The only matter needing attention 
is that we have enough topological slack in the deformation type to develop Gieseker 
bubbles at all $\pi$-splitting edges.

The weights are checked near $F_{\delta,z}$, where a formal neighborhood is isomorphic 
to the product of a scheme $J$ and several copies of the local models (Example \ref{Ex04}, 
Notation~\ref{localmodel} ).  The germ of $Z_\delta$ at $F_{\delta,z}$ 
consists of strata which are obtained by smooth Gieseker bubbles in a way which increases 
$d_+(\dt)$. It therefore lies over (the product of $J$ and) the $z$-axes of the local model, 
which tells us that the weights of $\Cs_+$ are positive.
\end{proof}

From the Bialynicki-Birula theorem \cite{BB}, we conclude the following. 
\begin{corollary}
The scalings $\eta_{0,\delta}: Z_\delta \to F_{\delta,z}$ and $\eta_{\infty,\delta}: 
W_\delta \to F_{\delta,w}$ are structure maps of bundles of affine spaces. 
In particular, $Z_\delta$ and $W_\delta$ are smooth.\qed
\end{corollary}

\section{Admissible Classes}\label{admissibleclasses}

We now estimate the weights of the admissible $K$-theory classes  over fixed points 
$F_n(\pi) \subset A_D$ as functions of the degrees $n_\pm$ (Remark~\ref{fixdegrees}).

\subsection{Definitions}
We recall the notation of Sec.~\ref{universal}: for any finite-dimensional 
representation $V$ of $\Cs$, $\phi^*V$ will be the vector bundle on $C$ 
associated to $V$ by $\mP$. Also recall the following (complexes of) coherent 
sheaves on $\Mtwid_{g,I}([\pt/\Cs])$: 
\begin{enumerate}
\item The {\it evaluation bundle} $\ev_i^*[V]=\sigma_i^*\phi^*V$.
\item Its {\it descendant bundles} are $\ev_i^*[V] \otimes [T_i^{\otimes j_i}]$, where 
$T_i = \sigma_i^*T_{\pi}$ is the relative tangent line to $C$ at $\sigma_i$, 
and $j_i$ is an integer.
\item The {\it Dolbeault\,\footnote{In \cite{math.AG/0312154}, the Dirac index class was used, 
but that requires a Spin structure on $C$.} index $I_V$ of $V$}, the 
complex $R\pi_*\phi^*V$. 
\item The {\it admissible line bundles} $\mL$:  negative 
(possibly fractional) powers of the determinant of cohomology 
of the standard representation $\C_1$ of $\Cs$,
\[
\mL \simeq (\det R\pi_*\phi^*\C_1)^{\otimes(-q)}, \quad q\in\Q_{>0}.
\]
 
\end{enumerate}

\begin{definition}
An {\it admissible complex} $\alpha$ on $\Mtwid_{g,I}([\pt/\Cs])$ is the tensor product 
of an admissible line bundle $\mL$ with any number of Dolbeault 
index and evaluation/descendant bundles.
\[
\alpha = \mL \bigotimes \otimes_a (R\pi_*\phi^*V_{a}) \bigotimes \otimes_b (\ev_b^*W_{b} 
\otimes T_i^{\otimes n_i}).
\]
{\it Admissible classes} are the topological K-theory classes of sums of admissible 
complexes.   
\end{definition}

\subsection{Weight Estimates in $A_D$}

Admissible complexes are bounded and coherent, so we can represent them locally 
as complexes $\mc{V}^\bullet$ of vector bundles equivariant under 
$\mc{G}(\pi)= \Cs_+\times\Cs_-$.  

\begin{proposition}
Fix $\pi \in \Pi(V)$. The following applies to the $\Cs_+\times\Cs_-$-weights of the 
fibers over $F_n(\pi)$, as $n_\pm$ vary:
\begin{enumerate}
\item For an evaluation or descendant class, they are bounded functions of $n_\pm$.
\item For an index complex $R\pi_*\phi^*V$, they are bounded functions of $n_\pm$,  
in a well chosen local resolution by vector bundles.
\item For an admissible line bundle $\mL$, they vary linearly with $n_\pm$, with 
positive coefficients.
\end{enumerate}
\end{proposition}

We  handle each case in turn. Let $C_f$ be the curve over $f \in F_n(\pi)$
and $\C_\lambda$ the irreducible $\Cs$-representation of weight $\lambda$. 

\begin{lemma}\label{lembundleweights}
  
Let $U$ be an open subset of an irreducible component $C'$ of $C_f$. 
\begin{enumerate}
\item If $C'$ is labeled by $\dt(\pi_+)$, then $\Cs_+$ acts on 
$\Gamma(U,\phi^*\C_\lambda)$ with weight $-\lambda$.  
\item If $C'$ is labeled by $\dt(\pi_-)$, then $\Cs_+$ acts on 
$\Gamma(U,\phi^*\C_\lambda)$ with weight $0$. 
\item Finally, if $U$ is a splitting Gieseker bubble, then both 
weights occur.
\end{enumerate}
\end{lemma}

\begin{proof}
Scaling all trivializations $t_v$, $v\in \pi_+$, by $g_+ \in \Cs$ can be 
absorbed by a global $g_+^{-1}$-rescaling the fibers of $\mc{P}$ on the $\dt(\pi_+)
$-components of $C_f$. Local sections of $\phi^*\C_\lambda$ on $U$ then rescale by 
$g_+^{-\lambda}$ or $1$, as appropriate. The Gieseker 
bubble is left to the reader (but see Remark~\ref{automph}). 
\end{proof}

\subsubsection{Evaluation and Descendant Classes}
At $f$, $\ev_i^*V$ is the fiber of $\phi^*V$ at $\sigma_i(f)$, and Lemma~\ref{lembundleweights} 
shows that the $\mc{G}(\pi)$-weights are independent of $n$. Further, $\mc{G}(\pi)$ 
acts trivially on the stable components of $C_f$, so the weights on the descendant 
class $\sigma_i^*(\phi^*V \otimes (T^*_\pi)^{\otimes j})$ are also constant in $n$.

\subsubsection{Index Classes}\label{indexweights}
A local complex of vector bundles $\mc{V}^\bullet$ representing $R^if_*\phi^*V$ 
can be built from sections $\{s_\alpha\}$ of a \v{C}ech resolution of $\phi^*V$. 
The fibers at $f$ of the $\mc{V}^i$ are spanned by the images of the 
generators $s_\alpha$: these are local sections of $\phi^*V|_{C_f}$, 
so Lemma~\ref{lembundleweights} above implies that their $\mc{G}(\pi)$-weights
don't depend on $n_\pm$.
\begin{remark}
If $\lambda\ge 0$, we need not refine the Gieseker bubbles to a \v{C}ech covering, 
because they give no $H^{1}$. Serre duality settles the case of  
$\lambda\le 0$. 
\end{remark}

\subsubsection{Admissible Line Bundles}

It suffices to compute the $\mc{G}(\pi)$-fixed point weights for $\det^{-1}R\pi_*\phi^*\C_1$. 
Flatness of $C$ implies that $\det^{-1}$ specializes at $f$ to
\[
\det\nolimits^{-1} R\Gamma(C_f;\phi^*\C_1). 
\]
We can exclude the splitting Gieseker bubbles in computing cohomology over $C_f$, 
and Lemma~\ref{lembundleweights} gives the $\mc{G}(\pi)$-character
\[
 (n_+ -g(\pi_+) +k(\pi_+))t_+^{-1} + (n_- -g(\pi_-)+k(\pi_-))t_-^{-1}.
\]
Thus, the $\Cs_\pm$-characters of $\det^{-q}$ are $t_\pm^{q(n_\pm -g(\pi_\pm) +k(\pi_\pm))}$.

\section{The Coherence Theorem}\label{theinvariants}
Here, we assemble the proof of our main theorem.
\begin{theorem}\label{maintheorem2}
The derived pushforward $R^\bullet F_*\alpha$ of any admissible complex $\alpha$ 
along $F:\Mtwid_{g,I}([\pt/\Cs]) \to \Mbar_{g,I}$  is coherent. 
\end{theorem}

\begin{proof}[Plan of proof]
We will check coherence in our \'etale charts $B$ of $\Mbar_{g,I}$.
Since $\Mtwid|_B\simeq A/\mc{G}$, coherence amounts to the local (over $B$) finite 
generation of $\mc{G}$-invariants in the derived global sections $R\Gamma(A,\alpha)$.

We prove this in two steps.  First, we fix the total degree $D$, and show that 
the $\mc{G}$-invariants in the derived global sections $R\Gamma(A_D,\alpha)$ are 
coherent.  Then we show that these invariants vanish for all but finitely many $D$.
\end{proof}

\subsection{Coherence on $A_D$}

If our base curve $\Sigma_o$ is reducible, $A_D$ is not proper over $B$: it has 
infinitely many finite-type strata.  We show that most of them do not contribute 
to the $\mc{G}$-invariants in $R\Gamma(A_D,\alpha)$. 

\begin{proposition}\label{vanishing}
Let $\mc{V}$ be a finite rank $\mc{G}$-equivariant vector bundle on $A_D$ with the 
following property:  

For all $\pi \in \Pi(V)$, the $\Cs_\pm$-weights of the 
fibers of $\mc{V}$ over the $\mc{G}(\pi)$-fixed points $F_n(\pi)$ are bounded below by 
increasing linear functions of $n_\pm$.  

Then, the $\mc{G}$-invariants in the local cohomology groups
\[
R^p\Gamma_{Z_\delta(\pi)}(A_D,\mc{V}) \qquad \mbox{and} 
	\qquad R^p\Gamma_{W_\delta(\pi)}(A_D,\mc{V})
\]
are finitely generated.  Moreover, the cohomologies vanish when $\delta\gg 0$.
\end{proposition}
\begin{proof}
We abbreviate $Z = Z_\delta(\pi)$ and $F = F_{\delta,z}(\pi)$ for fixed $\pi$;
the arguments for $W$ and $Z$ are similar, so we focus on $Z$. 

Now, $\mc{V}$ is a vector bundle and $Z$ is a smooth, closed subvariety of 
some open subspace $U\subset A_D$. Exactness of the functor of $\mc{G}$-invariants reduces the vanishing of the invariants in the cohomology groups with supports, $R^i\Gamma_{Z}(A_D,\mc{V})$, to that of the $\mc{G}$-invariants in $R^i\Gamma(U,R^\bullet\Gamma_{Z}(\mc{V}))$.  The latter will follow (via the filtration spectral sequence) from the vanishing of invariants in $R^i\Gamma(Z,\mc{V}\otimes \det N_{Z/A_D}\otimes \Sym N_{Z/A_D})$.

$Z$ is the total space of a bundle of affine spaces over the fixed point locus $F$, so by pushing down along the fibers we reduce computation of the latter to
\[
 R^i\Gamma(F, \mc{V}\otimes \det N_{Z/A_D}\otimes\Sym N_{Z/A_D} \otimes \Sym N^\vee_{F/Z}).
\]
The vector spaces in the two $\Sym$s in the RHS above have positive $\Cs_+$-weights.  Since $\mc{V}$ has finite rank, it follows that the $\mc{G}$-invariants in the RHS are finitely-generated.  Moreover, since $n_- \sim -\delta$ in $Z$, $n_{+} \gg 0$ if 
$\delta \gg 0$, so the $\Cs_{+}$-invariants vanish in that case.  Thus, the 
$\mc{G}(\pi)$-invariants in the RHS vanish, which implies that the $\mc{G}$-invariants 
vanish.
\end{proof}

The invariants in $R\Gamma(A^o_D,\alpha)$ are the direct images of $\alpha^{\Cs_\Delta}$ from $Q$, and the latter is proper over $B$ (Theorem~\ref{lbcoherence}); this ensures their 
coherence on $B$. Finite-generation of invariants in the local cohomologies, and their 
vanishing for almost all $\delta$, allows us to add arbitrarily many strata $Z_\delta(\pi)$ 
and $W_\delta(\pi)$ to $A^o_D$ without changing the finite-generation of invariants.

\begin{corollary}
The $\mc{G}$-invariants in $R\Gamma(A_D,\alpha)$ are finitely-generated.
\end{corollary}

\subsection{Varying $D$}
The diagonal subgroup $\Cs_\Delta \subset \mc{G}$ fixes every $A_D$, so its 
action on $R\Gamma(A_D,\alpha)$ comes from the fiber-wise action on $\alpha$.  
The $\Cs_\Delta$-weights appearing in this complex lie in a finite range; the width 
of this range depends on the class $\alpha$, but the upper and lower bounds  
grow linearly in $D$  ($q(n_++n_-) = qD + \mbox{constant}$), as enforced by the 
admissible line bundle factor in $\alpha$.
The $\mc{G}$-invariants in $R\Gamma(A_D,\alpha)$ therefore vanish for 
all but finitely many $D$.

\section{Towards Gromov-Witten Invariants for $[X/\Cs]$}\label{XmodC}

In this section, we assemble the technical ingredients that prepare the construction 
of Gromov-Witten invariants for quotient stacks $[X/\Cs]$. The material follows 
existing literature, with the modifications imposed by the $\Cs$-action.

\subsection{Definitions}
The yoga of stacks interprets a morphism  $\phi: Y\to [X/\Cs]$ from a scheme $Y$ to 
the quotient stack as a principal $\Cs$-bundle $\mP\to Y$, together with a section $s$ of 
the associated fiber bundle $\mP \times_{\Cs} X$. This certainly induces a continuous 
map of $Y$ into the Borel construction $X_{\Cs}$ of the quotient. If we set
\[
H_n([X/\Cs]) := H_{n}(X_{\Cs}),
\]
a map $\phi: C \to [X/\Cs]$ from an irreducible curve $C$ will have a definite {\it degree} 
$\beta \in H_2([X/\Cs])$. Reducible curves have a multi-degree, whose components 
sum up to the total degree.

\begin{definition}
A {\it Gieseker map from $C$ to $[X/\Cs]$} is a triplet $((C,\sigma_i),\mP,s)$ consisting of:
\begin{enumerate}
\item a prestable  marked curve $(C,\sigma_i)$,
\item a principal $\Cs$-bundle $p: \mP \to C$, and
\item A $\Cs$-equivariant map $s: \mP \to X$ (equivalently, a section of the associated bundle $\XP = \mP \times_{\Cs} X$ with fiber $X$)
\end{enumerate}
such that
\begin{enumerate}
\item $\mP$ has degree $0$ on any irreducible rational component of $C$ which has one 
node and one marked point.
\item $\mP$ has either degree $0$ or degree $1$ on any unstable rational component 
of $C$ which has two nodes.
\item $s$ is non-trivial on any unstable component on which $\mP$ has degree $0$.
\end{enumerate}

We denote by $\Mtwid_{g,I}([X/\Cs])$ the fibered category of Gieseker maps to $[X/\Cs]$ 
from stable marked curves of type $(g,I)$.  Its connected components carry definite total 
degrees $\beta \in H_2([X/\Cs])$, and we assemble maps of degree $\beta$ into 
$\Mtwid_{g,I,\beta}([X/\Cs])$.
\end{definition}

There is a forgetful map
\[
F_s: \Mtwid_{g,I,\beta}([X/\Cs]) \to \Mtwid_{g,I,ft_*\beta}(\pt/\Cs),
\]
where $ft_*\beta$ is the degree obtained from the homomorphism $ft_*: H_2([X/\Cs]) 
\to H_2([\pt/\Cs])$ (where we should not forget to contract any rational components 
which carry trivial bundles, after forgetting the section).

\subsection{Properness}

\begin{theorem}\label{propdm}
$F_s$ is proper and Deligne-Mumford.
\end{theorem}

\begin{proof}
In essence, this is because $[X/\Cs] \to [\pt/\Cs]$ is proper and representable.  
We make the argument precise via the valuative criteria for completeness \& separability.  

Suppose that we have a family $(C,\sigma_i,\mP,s)$ over the punctured disk $D^\times$ and 
an extension $(Z,z_i,\mc{R})$ of $F_s(C,\sigma_i,\mP)$ to the disk $D$. 
(Any required base change will be subsumed in the notation.)

\emph{Completeness}:  We want an extension $(Y,y_i,\mc{Q},t)$ of $(C,\sigma_i,\mP,s)$ 
to $D$ such that  $F_s(Y,y_i,\mc{Q},t) = (Z,z_i,\mc{R})$.

First, we extend the family $(C,\sigma_i)$ to $D$.  This may require base change, and is an easy consequence of the existence of nodal reduction \cite{MR1631825}.  We denote the extension by $(Y',y_i')$; it comes equipped with a contraction map $c: Y' \to Z$.   We denote by $\mc{Q}'$ the pullback $c^*\mc{R}$; note that $\mc{Q}'$ is trivial on components collapses by $c$.

Given $Y'$, the graph of $s$ gives us an embedding $j: C \to X_{\mc{Q}'}$, where $X_{\mc{Q}'}$ is the associated bundle $\mc{Q}' \times_{\Cs} X$.  The morphism $u: X_{\mc{Q}'} \to B$ has compact fibers, so the closure $\overline{j(C)}$ of the image of $j$ is also a finite type curve over $B$.  $\overline{j(C)}$ is not necessarily prestable.  However, resolution of singularities leads to a prestable curve $Y''$ (with a 
resolution map $r: Y'' \to \overline{j(C)}$); base change may also be required 
at this step.  This gives us a sequence of maps (over $D$)
\[
\xymatrix{Y'' \ar[r]^r & \overline{j(C)} \ar@{^{(}->}[r]^j & X_{\mc{Q}'} \ar[r]^{pr} & Y' \ar[r]^{c} & Z},
\]
where $pr: X_{c^*\mc{R}} \to \Sigma'_0$ is the bundle structure map.  The composition 
$c_r=c\circ pr \circ j \circ r: Y'' \to Z$ is a contraction map. We denote the pullback $c_r^*\mc{R}$ by $\mc{Q}''$, and the lifts of the marked points $z_i$ by $y_i''$.

Pulling back $\mc{R}$ step by step from $C$ to $Y''$, we get a sequence of bundles, the last of which is $c_r^*\mP$, as in the diagram below.
\[
\xymatrix{
\mc{Q}'' \ar[rr] \ar[d] &   & \mc{Q}' \times X \ar[r]^{pr_1} \ar[d] &\mc{Q}' \ar[d] \ar[r] & \mP \ar[d]\\
Y'' \ar[r]^r & \overline{j(C)} \ar@{^{(}->}[r]^j & X_{\mc{Q}'} \ar[r]^{pr} & Y' \ar[r]^{c} & Z}
\]
We also get a section $s': \mc{Q}''\to X$ from the composition 
\[
\mc{Q}'' \to \mc{Q}' \times X \to X.
\]  
The collection $(Y'',y_i'',\mc{Q}'',s')$ is a map to $[X/\Cs]$, but not necessarily a Gieseker map, as the curve may have unstable components carrying a trivial bundle and a trivial section.  We obtain the desired extension by contracting these unstable components.

\emph {Separability}:  
Suppose that we are given two different pairs $(Y_1,y_{1,i},\mc{Q}_1,s_1)$ and $(Y_2,y_{2,i},\mc{Q}_2,s_2)$ which both extend the given family over $B$ compatibly with the given Gieseker map $(Z,z_i,\mc{R})$ to $[\pt/\Cs]$.  We may freely suppose that both extensions are defined over the same base extension.

Consider the fiber product $Y_1 \times_Z Y_2$. Our assumptions imply that $Y_1|_{B^\times} = Y_2|_{B^\times}$ and that the special fibers of $Y_1$ and $Y_2$ both contract onto the special fiber of $Z$.  It follows that all the maps in the bottom diamond of the diagram below are contraction maps.
\begin{equation*}
\xymatrix{
& \mc{Q} \ar[dr]^{f_2} \ar[dl]_{f} \ar[d] & \\
\mc{Q}_1 \ar[d] \ar[dr]& Y_1 \times_{Z} Y_2\ar[dr] \ar[dl] & \mc{Q}_2 \ar[d] \ar[dl]\\
Y_1 \ar[dr] & \mc{R}\ar[d]  & Y_2 \ar[dl]\\
 & Z& }
\end{equation*}
Moreover the two sections $\mc{Q} \to \mc{Q}_1 \to X$ and $\mc{Q} \to \mc{Q}_2 \to X$ agree on the open dense set $\mc{Q}|_{B^\times}$.  $X$ is separated, so the two sections agree.  The Gieseker map obtained by contracting any unstable components in $Y_1\times_{Z}Y_2$ is unique, so it follows that the two given families are isomorphic.

\emph {Deligne-Mumford}:  Let $C_v$ be a component of $C$.  If $C_v$ is contracted to a point by the section-forgetting morphism $F_s$, then $\mP|_{C_v}$ is trivial, so $s|_{C_v}$ must be equivalent to a non-trivial map $C_v \to X$.  We know from Gromov-Witten theory that such maps admit only finitely many automorphisms.
 
On the other hand, if $C_v$ is stable, then the existence of a non-trivial section on $C$ can only reduce the number of automorphisms. 
\end{proof}

\subsection{Virtual smoothness}
Forgetting the section $s$, but not contracting unstable components defines a morphism to 
the stack $\mf{M}_{g,I,ft_*\beta}([\pt/\Cs])$ of bundles on all prestable curves of degree $ft_*\beta$:  
\[
\widetilde{F}_s: \Mtwid_{g,I,\beta}([X/\Cs]) \to \mf{M}_{g,I,ft_*\beta}([\pt/\Cs]).
\]

\begin{theorem}\label{perfobs}
$L_{\widetilde{F_s}}$ admits a relative perfect obstruction theory.
\end{theorem}
\noindent
Recall from \cite{MR1437495} that a {\it relative perfect obstruction theory} for 
the cotangent complex $L_{\widetilde{F_s}}$ is pair $(E,e)$ consisting of an element 
$E$ of the derived category of $\Mtwid_{g,I,\beta}([X/\Cs])$, and a homomorphism 
$e: E \to L_{\widetilde{F}_s}$ in the derived category, such that
\begin{enumerate}
\item $E = [E^{-1} \to E^0]$ is locally equivalent to a two-term complex
  of locally free sheaves.
\item $H^0(e)$ is an isomorphism.
\item $H^{-1}(e)$ is a surjection.
\end{enumerate}

\begin{proof}[Proof of Theorem \ref{perfobs}]
The proof is an almost word-for-word copy of the one given by 
Behrend \& Fantecchi in \cite{MR1431140, MR1437495}.

Fix a curve $C$ and a principal $\Cs$-bundle $p:\mP \to C$,
and let $\Gamma$ denote the space $\op{Hom}_{\Cs}(\mP,X)$ of sections.
$\Gamma$ comes equipped with universal families:
\begin{equation*}
\xymatrix{ 
\mP\times \Gamma \ar[r]^{s} \ar[d]^{p \times \op{id}_\Gamma} & X \ar[d]^{\rho}\\
C \times \Gamma \ar[r]^{\phi_s} \ar[d]^{\pi} & [X/\Cs]\\
\Gamma & }
\end{equation*}
It follows from the functorial properties of the cotangent complex that
we have a morphism $\tilde{e}: s^*L_X \to p^*\pi^*L_\Gamma$.   If we
take $\Cs$-invariants in the pushdown via $p$, we get
\begin{equation*}
\tilde{e}': (p_*s^*L_X)^{\Cs} \to \pi^*L_\Gamma.
\end{equation*}
Tensoring with the dualizing complex of $C$, we obtain a morphism 
\[
\tilde{e}'': \omega_{C}\otimes (p_*s^*L_X)^{\Cs} \to
\omega_{C}\otimes \pi^*L_\Gamma = \pi^!L_\Gamma.
\]
Then, by adjunction, we have a morphism 
\[
\tilde{e}'': R\pi_*(\omega_{C}\otimes (p_*s^*L_X)^{\Cs}) \to L_\Gamma.
\]
Finally, it follows from Verdier duality that
\[
R\pi_{*}(\omega_{C}\otimes (p_*s^*L_X)^{\Cs}) =R\pi_*(p_*s^*T_X)^{\Cs},
\] 
and so we have a morphism
\[
e: [R\pi_*(p_*s^*T_{X})^{\Cs}]^{\vee} \to L_\Gamma.
\]
This morphism is a perfect obstruction theory for $L_\Gamma$; the proof
is as in \cite{MR1437495}.  Moreover, all of the
objects here generalize well to the relative case, and therefore apply
to the universal family.  Thus, we have a perfect relative obstruction
theory
\[
e: E = [R\pi_*(p_*s^*T_{X})^{\Cs}]^{\vee} \to L_{\widetilde{F_s}},
\]
where now $\pi$, $p$, and $s$ refer to the universal families on the
moduli stack.
\end{proof}

Given this perfect obstruction theory, the virtual normal cone device developed by 
Behrend \& Fantecchi \cite{MR1437495} constructs the {\it virtual structure sheaf} 
$\mO^{vir}$ in the bounded derived category of coherent sheaves on $\Mtwid_{g,I,\beta}([X/\Cs])$. 
This represents a family of virtual fundamental K-homology cycles, which we use to integrate 
along the fibers of $F_s$. A construction closely suited to our purposes is found in \cite{MR2040281}.

For a vector bundle $\mc{V}$ on $\Mtwid_{g,I,\beta}([X/\Cs])$, let $\mc{V}^{vir}: = 
\mc{V}\otimes\mc{O}^{vir}$ and define the {\it virtual pushforward of $\mc{V}$} 
to $\Mtwid_{g,I}(\pt/\Cs)$ as
\[
(F_s)^{vir}_![\mc{V}] = (F_s)_*[\mc{V}^{vir}].
\]
Thanks to Theorem~\ref{propdm}, this has a well-defined $K$-theory class.  

\subsection{Gromov-Witten invariants?}
There is an obvious notion of admissible class on $\Mtwid_{g,I,\beta}([X/\Cs])$, produced 
from tautological classes and line bundle twists, whose direct images should give Gromov-Witten 
invariants for $[X/\Cs]$. This requires a finiteness for the virtual pushforward analogous 
to the main theorem of this paper. While we do not know if the $(F_s)_*^{vir}$ above  
takes such classes to admissible classes on $\Mtwid_{g,I}([\pt/\Cs])$, 
we expect that the result satisfies the weight bounds of Lemma~\ref{vanishing}. The same 
argument would then ensure finiteness.

\def\cprime{$'$} \def\cprime{$'$}

\end{document}